\documentclass[10pt]{article}
\date{{\tt diffvect3.tex}, November 30, 2010} 
\usepackage{amssymb}
\usepackage{amsmath}
\input{liemacs.sty} 

\begin{document} 



\title{On Differentiable Vectors for \\ 
Representations of Infinite Dimensional Lie Groups} 
\author{Karl-Hermann Neeb\begin{footnote}{
Department  Mathematik, FAU Erlangen-N\"urnberg, Bismarckstrasse 1 1/2, 
91054-Erlangen, Germany; neeb@mi.uni-erlangen.de}
\end{footnote}
\begin{footnote}{Supported by DFG-grant NE 413/7-1, Schwerpunktprogramm 
``Darstellungstheorie''.} 
\end{footnote}}

\maketitle
\date

\begin{abstract} 
In this paper we develop two types of tools to deal with differentiability 
properties of vectors in continuous representations $\pi \: G \to \GL(V)$ 
of an infinite dimensional Lie group $G$ on a locally convex space~$V$. 
The first class of results concerns the space 
$V^\infty$ of smooth vectors. If $G$ is a Banach--Lie group, 
we define a topology on the space 
$V^\infty$ of smooth vectors for which the action of $G$ on this 
space is smooth. If $V$ is a Banach space, then 
$V^\infty$ is a Fr\'echet space. This applies in particular 
to $C^*$-dynamical systems $(\cA,G, \alpha)$, where 
$G$ is a Banach--Lie group. 
For unitary representations we show that 
a vector $v$ is smooth if the corresponding 
positive definite function $\la \pi(g)v,v\ra$ is smooth. 

The second class of results concerns criteria for 
$C^k$-vectors in terms of operators of the derived representation 
for a Banach--Lie group $G$ acting on a Banach space~$V$. 
In particular, we provide for each $k \in \N$ 
examples of continuous unitary representations  
for which the space of $C^{k+1}$-vectors is trivial and the space 
of $C^k$-vectors is dense. 
\end{abstract} 

\tableofcontents

\section{Introduction} \label{sec:1}

Let $G$ be a Lie group modeled on a 
locally convex space (cf.~\cite{Ne06} for a survey on locally convex 
Lie theory). A representation $\pi \: G \to \GL(V)$ of $G$ on 
the locally convex space $V$ is called {\it continuous} if 
it defines a continuous action of $G$ on $V$. We call 
an element $v \in V$ a {\it $C^k$-vector} ($k \in \N_0 \cup\{\infty\}$) 
if the orbit map $\pi^v \: G \to V, g \mapsto \pi(g)v$, is a 
$C^k$-map and $V^k = V^k(\pi)$ for the space of $C^k$-vectors 
(cf.\ Section~\ref{sec:2} for precise definitions). 

It is a fundamental problem in the representation theory of Lie groups 
to understand the interplay between differentiability 
properties of group representations and representations of 
the Lie algebra $\g= \L(G)$ of $G$. In particular, the following 
questions are of interest: 
\begin{description}
\item[\rm(1)] How to find $C^k$-vectors? Can they be characterized 
in terms of the Lie algebra? 
\item[\rm(2)] Does the space $V^\infty$ of smooth vectors 
carry a natural topology for which the action of $G$ on 
this space is smooth? 
\item[\rm(3)] Under which circumstances do differentiable, resp., smooth 
vectors exist, resp., form a dense subspace? 
\end{description}

In this note we provide answers to some of these questions 
for infinite dimensional Lie groups~$G$. For finite dimensional 
groups most of these results are either trivial or well-known 
(see \cite{Go69} for (1), \cite{Mi90} for (2) and \cite{Ga47} for 
(3)), so that 
the main issue is to identify and deal with the subtleties of infinite 
dimensional Lie groups.

After discussing some preliminaries on Lie groups and their representations 
in Sections~\ref{sec:2} and \ref{sec:3}, we turn in 
Section~\ref{sec:4} to the space $V^\infty$ of smooth vectors of the 
representation of a Banach--Lie group $G$ on a locally convex 
space~$V$. Here our main result is that, by embedding $V^\infty$ 
into a product of the Banach spaces $\Mult^p(\g,V)$ of 
continuous $p$-linear maps $\g^p\to V$, we obtain a locally convex 
topology on $V^\infty$ for which the $G$-action on $V^\infty$ 
is smooth (Theorem~\ref{thm:4.1.3}). For finite dimensional groups, 
this follows easily from the smoothness of the translation 
action of $G$ on the space $C^\infty(G,V)$ (cf.\ Proposition~\ref{prop:findim}), 
but for infinite dimensional groups there seems to be no ``good'' topology 
on the space $C^\infty(G,V)$. The smooth compact open topology is too 
coarse to ensure continuity of the action. For 
the Fr\'echet--Lie group $G = \R^\N$ and the unitary representation 
on $\cH = \ell^2(\N,\C)$ defined by $\big(\pi(x)y\big)_n = e^{ix_n}y_n$,  
there exists no locally convex topology on $\cH^\infty  = \C^{(\N)}$ 
for which the $G$-action is smooth. In view of this example, 
Theorem~\ref{thm:4.1.3} 
on the smoothness of the action on $V^\infty$ does not 
extend to unitary representations of Fr\'echet--Lie groups. 
The locally convex space  $V^\infty$ has a topological dual space 
$V^{-\infty} := (V^\infty)'$, whose elements are called 
{\it distribution vectors}. For finite dimensional groups they 
correspond to equivariant embeddings of $V$ into the space of 
distributions on $G$, and in our context they still parametrize 
equivariant embeddings $V^\infty \to C^\infty(G)$, where $C^\infty(G)$ 
carries the left regular representation. In this sense 
the spaces $V^{-\infty}$ provide natural candidates for spaces 
of distributions on Banach--Lie groups. For discussions of 
distributions on infinite dimensional vector spaces, 
we refer to \cite{KSWY02}, \cite{KS76} and \cite{Ko80}. 

We continue the discussion of $V^\infty$ in Section~\ref{sec:5}, 
where we show that if $V$ is a Banach space, 
then $V^\infty$ is complete, hence a Fr\'echet space. 
In Section~\ref{sec:6} we briefly discuss an 
application of these results to $C^*$-dynamical systems 
$(\cA,G,\alpha)$, where $G$ is a Banach--Lie group and 
$\cA$ is a unital Banach algebra. Then the Fr\'echet space 
$\cA^\infty$ is shown to be a continuous inverse algebra, i.e., 
its unit group is open and the inversion is a continuous map. 
The main result of Section~\ref{sec:7} is that, for a unitary 
representation $(\pi, \cH)$, a vector $v \in \cH$ is smooth if and only 
if the corresponding matrix coefficient $\pi^{v,v}(g) = \la \pi(g)v,v\ra$ 
is smooth. 

Although smooth vectors form the natural domain for the derived 
representation of the Lie algebra $\g$ of $G$, in many situations it is 
desirable to have some information on 
$C^k$-vectors for $k < \infty$. This motivates the discussion of 
$C^1$-vectors for continuous representations of 
Banach--Lie groups on Banach spaces 
in Section~\ref{sec:9}. 
The main result is that the intersection 
$\cD_\g \subeq V$ of the domains of the generators 
$\oline{\dd\pi}(x)$ of the one-parameter groups 
$t \mapsto \pi(\exp_G(tx))$ on $V$ coincides with the space 
of $C^1$-vectors (Theorem~\ref{thm:contlin}). Here the main difficulty 
lies in the verification of the additivity of the map 
$\omega_v \: \g \to V, x \mapsto \oline{\dd\pi}(x)v$, 
for which we use quite 
recent refinements of Chernoff's Theorem by 
Neklyudov (\cite{Nek08}). If $\omega_v$ is assumed to be continuous, 
its linearity is easily verified (Lemma~\ref{lem:a.11c}). 
It is also easily verified for finite dimensional groups (\cite[p.~221]{Go70}). 
All this is applied in Section~\ref{sec:10} 
to $C^k$-vectors, which are shown to coincide (for a representation 
of a Banach--Lie group on a Banach space) with the space 
$\cD_\g^k$, the common domain of $k$-fold products of the operators 
$\oline{\dd\pi}(x)$, $x \in \g$ (Theorem~\ref{thm:10}). 
Here an interesting point is that the $C^k$-concept we use is weaker
than the $C^k$-concept used in the context of Fr\'echet differentiable 
maps between Banach spaces. The discussion of the examples 
in Section~\ref{sec:11} shows that, in general, the space of 
Fr\'echet-$C^k$-vectors is a proper subspace of $\cD_\g^k$ which can 
be trivial. 

A crucial problem one has to face in the representation theory 
of infinite dimensional Lie groups is that a continuous 
representation need not have any differentiable vector, 
as is the case for finite dimensional Lie groups (\cite{Ga47}). 
Refining a construction from \cite{BN08}, we discuss in 
Section~\ref{sec:11} the continuous unitary representation 
of the Banach--Lie group $G = (L^p([0,1],\R),+)$ 
on the Hilbert  space $L^2([0,1],\C)$ by $\pi(g)f = e^{ig}f$ 
and determine for every $k \in \N$ the space of 
$C^k$-vectors. For $p = 1$ we thus obtain a continuous 
representation with no non-zero $C^1$-vector, and for 
$p = 2k$, the space of $C^k$-vectors is non-zero, but 
there is no non-zero $C^{k+1}$-vector.

The situation is much better for Lie groups which are 
direct limits of finite dimensional ones. For these 
groups strong existence results on smooth vectors are available 
(\cite{Sa91}, \cite{Da96} for unitary representations and 
\cite{Si75} for Banach representations), 
and we explain in Section~\ref{sec:12} how they fit into 
our general framework. 
The density of $\cH^\infty$ 
for particular unitary representations of diffeomorphism groups 
is shown in \cite{Sh01}. 
In Section~\ref{sec:13} we show that, for rather trivial 
reasons, the space of smooth vectors of a continuous unitary 
representation is also dense for projective limits 
of finite dimensional Lie groups. 
We conclude this paper in Section~\ref{sec:14} with a discussion 
of some open problems. 

The class of groups for which the theory 
developed in this article has the strongest impact is the class
of Banach--Lie groups and we hope in particular that 
our results on smooth vectors will lead to a better understanding 
of their unitary representations. To  make this more concrete, 
we recall some of the major classes of Banach--Lie groups and 
what is known on their unitary representation theory. 

For a Banach--Lie group $G$, the most regular class of unitary 
representations are the {\it bounded} ones, i.e., those 
for which $\pi \: G \to \U(\cH)$ is a morphism of Banach--Lie groups 
with respect to the norm topology on $\U(\cH)$. It is easy to see 
that faithful representations of this kind exist only when 
$\L(G)$ carries an $\Ad(G)$-invariant norm, a property which for finite 
dimensional Lie algebras is equivalent to compactness. 
If $G$ has a universal complexification 
$\eta_G \: G \to G_\C$ (see \cite{GN03} for criteria on 
the existence), then every bounded unitary representation 
extends to a holomorphic representation 
$\hat\pi \: G_\C \to \GL(\cH)$. Holomorphic representations 
of Banach--Lie groups have been studied systematically 
in the context of the large class of root graded Lie groups 
in \cite{MNS10} and \cite{NS10}. The class of root graded Lie groups 
contains in particular analogs of classical groups over 
unital Banach algebras and the most typical example 
is $\GL_n(\cA)$ for a unital Banach algebra $\cA$. 
In \cite{MNS10} we developed a Borel--Weil theory for such 
groups, showing in particular that representations 
in spaces of holomorphic sections of bundles over 
the natural analogs $G/P$ of flag manifolds are always 
holomorphic representation by bounded operators and that all 
holomorphic irreducible representations are of this form. 
In this context the most intriguing problem is to decide 
for which holomorphic representations 
$\rho \: P \to \GL(V)$, the associated holomorphic vector 
bundle has nonzero holomorphic sections. For 
root graded groups whose Lie algebra has the form 
$\g \otimes \cA$, $\g$ complex semisimple and $\cA$ a 
unital commutative Banach algebra, this question has been 
answered completely in \cite{NS10} for the special case 
of line bundles, i.e., $\rho \: P \to \C^\times$ is a holomorphic 
character. From this classification result it follows 
in particular that the spaces of holomorphic sections of these 
line bundles are always finite dimensional. As a special 
case, this theory leads to a complete classification of the 
bounded unitary representations of the unitary groups 
$\U_n(\cA)$, where $\cA\cong C(X)$ is a commutative unital $C^*$-algebra 
and all these representation factor through some quotient group $\U_n(\C)^d$. 

In \cite{BN10} the Schur--Weyl theory of representations of 
$\U_n(\C)$ is extended to unitary groups $\U(\cA)$ of a 
$C^*$-algebra. It is shown that tensor products of 
irreducible representations of $\cA$ and their duals 
decompose into finitely many bounded unitary representations of 
$\U(\cA)$ according to the classical Schur--Weyl pattern. 
If $\cA$ is commutative, then the aforementioned results from 
\cite{NS10} imply that we thus obtain all irreducible unitary 
representations of $\SU_n(\cA) \subeq \U(M_n(\cA))$, 
but if $\cA$ is non-commutative the classification problem for 
irreducible unitary representations of $\U_n(\cA)$ is open. 
Thanks to the existence of group complexifications, 
the main tool in the analysis of bounded unitary representations of 
Banach--Lie groups is complex analysis, resp., holomorphic 
extension. This method is also exploited extensively in the classification 
of bounded unitary representations of the operators groups 
\[ \U_p(\cH) := \U(\cH) \cap (\1 + \cL_p(\cH)),\quad 1 \leq  p \leq \infty,\] 
where $\cL_p(\cH)$ is the {\it $p$-th Schatten ideal} 
(see also \cite{Ne04}). Here a remarkable result is that the 
classification does not depend on $p$ in the range $1 < p \leq \infty$, 
so that the picture is the same as for 
$p = \infty$, where $\cL_\infty(\cH) = K(\cH)$ is the 
$C^*$-algebra of compact operators, and the irreducible 
representations are obtained from Schur--Weyl theory by decomposing 
tensor products of the form $\cH^{\otimes n} \otimes \oline\cH^{\otimes m}$. 
For $p = 1$, the bounded representation 
theory of $\U_1(\cH)$ is much richer but not of type $I$. 
A classification of the bounded irreducible, resp., factor representations  
of this group is still an interesting open problem. 

In this context a remarkable result of D.~Pickrell asserts that 
for the full unitary group $\U(\cH)$ of an infinite dimensional separable Hilbert space 
over $\R, \C$ or $\bH$, every separable continuous unitary 
representation is automatically bounded, a direct sum of 
irreducible representations, and these are obtained 
by Schur--Weyl theory from the decomposition of the tensor 
products $\cH_\C^{\otimes n}$ (\cite{Pi88}). Basically this means that, 
as far as separable representations are concerned, the full 
unitary groups and finite products thereof behave exactly like 
compact groups. 

It is a rule of thumb that complex analytic methods apply 
well to bounded and semibounded unitary representations, 
where the latter class is defined by the condition that the operators 
$i\dd\pi(x)$ are uniformly bounded above on some open subset of the 
Lie algebra (see \cite{Ne10b} for a survey on semibounded representations). 
Beyond semibounded representations one has to rely on real analytic techniques. 

Typical analogs of finite dimensional non-compact semisimple 
groups are the automorphism groups of Hilbert--Riemannian 
symmetric spaces such as restricted Gra\3mannians and symmetric 
Hilbert domains (\cite{NO98}, \cite{Ne01c}, \cite{Ne04}).
For the same reason that finite dimensional 
groups with non-compact Lie algebra have no faithful 
finite dimensional unitary representation, 
these groups have no faithful bounded unitary representation. 
An important example is the {\it symplectic 
group} $\Sp(\cH)$ of all real linear symplectic automorphisms 
of $(\cH,\omega)$, where $\cH$ is a complex Hilbert space and 
$\omega(v,w) = \Im\la v,w\ra$. Even more important is the 
{\it restricted symplectic group} 
$\Sp_{\rm res}(\cH)$ consisting of all elements for which 
$g^\top g-\1$ is Hilbert--Schmidt (cf.\ \cite{Ne10b}). 
The latter group has a by far richer (projective) representation theory 
than the former. It is the prototype of a hermitian 
Lie groups, i.e., an automorphism group of an infinite 
dimensional hermitian symmetric space. 
For hermitian Lie groups 
and their central extensions, the irreducible semibounded unitary 
representations have recently been classified in \cite{Ne10c} 
by a combination of complex analytic methods with Pickrell's Theorem (\cite{Pi88}) and some results on unitarity of highest weight modules (\cite{NO98}). 
For $\Sp_{\rm res}(\cH)$ it turns out that all separable semibounded 
representations are trivial, but it has a central extension with a 
rich semibounded representation theory. For the full symplectic 
group $\Sp(\cH)$ it seems quite likely that all its unitary 
representation are trivial. So it becomes an interesting issue 
which structural features of a Banach--Lie group 
lead to obstructions for the existence of non-trivial unitary 
representations. 

An area which is still largely unexplored is the theory 
of $(\g,K)$-modules for pairs such as $(\sp_{\rm res}(\cH), \U(\cH))$. 
Related problems have been addressed by 
G.~Olshanski from the point of view of topological groups 
in terms of $(G,K)$-pairs, where $K$ is a subgroup of $G$, 
such as $(\Sp_{\rm res}(\cH), \U(\cH))$ (cf.\ \cite{Ol90}, \cite{Pi90}),  
and we hope that, eventually, a better understanding of 
differentiability properties as 
discussed in the present paper and integration techniques 
as in \cite{Ne10a} and \cite{Me10} will lead 
to a more transparent Lie theoretic 
understanding of the representations of such groups. 

The groups $G = C^k(M,K)$ of $C^k$-maps, $k \geq 1$, 
on a compact manifold $M$ with values in a finite dimensional Lie group~$K$ 
form a another class of interesting Banach--Lie groups; more generally 
one considers $C^k$-gauge transformations of principal bundles. 
These groups have interesting central extensions $\hat G$ 
(\cite{Ne10d}), and for  
the case $M = \bS^1$, where $G$ is a loop group, the central 
extension $\hat G$ has an interesting family of unitary representations 
extending to a semidirect product $\hat G \rtimes \T$, where $\T$ acts 
on $\bS^1$ by rigid rotations (cf.\ \cite{PS86}, \cite{Ne01b}, \cite{Alb93}). 
These groups are Banach manifolds 
and topological groups, but their multiplication is not smooth. 
We expect that a thorough analysis of such semidirect products 
$N \rtimes H$, where $H$ is a Banach--Lie group acting 
{\sl continuously} on the Banach--Lie group $N$, can be based 
on the results of the present paper. For the case where 
$(\cA, H,\alpha)$ is a $C^*$-dynamical system, the corresponding 
covariant representations lead in particular to unitary representations 
of the groups $\U(\cA)\rtimes_\alpha H$ which are of this type 
(cf.~\cite{Pe89}). 
Another interesting class of such groups providing testing cases 
for a general theory are the semidirect products 
$\Heis(\cH) \rtimes_\alpha \R$, where 
$\cH$ is a complex Hilbert space and $\alpha(t) = e^{itH}$ a 
strongly continuous unitary one-parameter group (cf. \cite{NZ10}). 

\section{Locally convex Lie groups} \label{sec:2}

In this section we briefly recall the basic concepts related 
to infinite dimensional Lie groups. 

\begin{defn} 
Let $E$ and $F$ be locally convex spaces, $U
\subeq E$ open and $f \: U \to F$ a map. Then the {\it derivative
  of $f$ at $x$ in the direction $h$} is defined as 
$$ \dd f(x)(h) := (\partial_h f)(x) := \derat0 f(x + t h) 
= \lim_{t \to 0} \frac{1}{t}(f(x+th) -f(x)) $$
whenever it exists. The function $f$ is called {\it differentiable at
  $x$} if $\dd f(x)(h)$ exists for all $h \in E$. It is called {\it
  continuously differentiable}, if it is differentiable at all
points of $U$ and 
$$ \dd f \: U \times E \to F, \quad (x,h) \mapsto \dd f(x)(h) $$
is a continuous map. Note that this implies that the maps 
$\dd f(x)$ are linear (cf.\ \cite[Lemma~2.2.14]{GN10}). 
The map $f$ is called a {\it $C^k$-map}, $k \in \N \cup \{\infty\}$, 
if it is continuous, the iterated directional derivatives 
$$ \dd^{j}f(x)(h_1,\ldots, h_j)
:= (\partial_{h_j} \cdots \partial_{h_1}f)(x) $$
exist for all integers $j \leq k$, $x \in U$ and $h_1,\ldots, h_j \in E$, 
and all maps $\dd^j f \: U \times E^j \to F$ are continuous. 
As usual, $C^\infty$-maps are called {\it smooth}. 
\end{defn}

Once the concept of a smooth function 
between open subsets of locally convex spaces is established  
(cf.\ \cite{Ne06}, \cite{Mil84}, \cite{GN10}), it is clear how to define 
a locally convex smooth manifold. 
A {\it (locally convex) Lie group} $G$ is a group equipped with a 
smooth manifold structure modeled on a locally convex space 
for which the group multiplication and the 
inversion are smooth maps. We write $\1 \in G$ for the identity element and 
$\lambda_g(x) = gx$, resp., $\rho_g(x) = xg$ for the left, resp.,
right multiplication on $G$. Then each $x \in T_\1(G)$ corresponds to
a unique left invariant vector field $x_l$ with 
$x_l(g) := T_\1(\lambda_g)x, g \in G.$
The space of left invariant vector fields is closed under the Lie
bracket of vector fields, hence inherits a Lie algebra structure. 
In this sense we obtain on $\g := T_\1(G)$ a continuous Lie bracket which
is uniquely determined by $[x,y]_l = [x_l, y_l]$ for $x,y \in \g$. 
We shall also use the functorial 
notation $\L(G) := (\g,[\cdot,\cdot])$ 
for the Lie algebra of $G$ and, accordingly, 
$\L(\phi) = T_\1(\phi)\: \L(G_1) \to \L(G_2)$ 
for the Lie algebra morphism associated to 
a morphism $\phi \: G_1 \to G_2$ of Lie groups. 
Then $\L$ defines a functor from the category 
of locally convex Lie groups to the category of locally convex topological 
Lie algebras. The adjoint action of $G$ on $\g$ is defined by 
$\Ad(g) := \L(c_g)$, where $c_g(x) = gxg^{-1}$ is the conjugation map. 
The adjoint action is smooth and 
each $\Ad(g)$ is a topological isomorphism of $\g$. 
If $\g$ is a Fr\'echet, resp., a Banach space, then 
$G$ is called a {\it Fr\'echet-}, resp., a 
{\it Banach--Lie group}. 

For every Lie group $G$, the tangent bundle $TG$ is a Lie group 
with respect to the tangent map $T(m_G)$ of the multiplication 
$m_G \: G\times G \to G$ on $G$. It contains $G$ as the zero section, 
which is a Lie subgroup, and the projection $TG \to G$ is a morphism of 
Lie groups whose kernel is the additive group 
of $\g \cong T_\1(G)$. In this sense we write 
$g.x = T_{(g,\1)}(m_G)(0,x) = T_\1(\lambda_g)x$ and 
$x.g = T_{(\1,g)}(m_G)(x,0) = T_\1(\rho_g)x$ for 
$g \in G$ and $x \in \g$. The two maps 
\begin{equation}
  \label{eq:tang-triv}
G \times \g \to TG, \quad (g,x) \mapsto g.x 
\quad \mbox{ and } \quad 
G \times \g \to TG, \quad (g,x) \mapsto x.g 
\end{equation}
trivialize the tangent bundle. 

A smooth map $\exp_G \: \g \to G$  is called an {\it exponential function} 
if each curve $\gamma_x(t) := \exp_G(tx)$ is a one-parameter group 
with $\gamma_x'(0)= x$. The Lie group $G$ is said to be 
{\it locally exponential} 
if it has an exponential function for which there is an open $0$-neighborhood 
$U$ in $\g$ mapped diffeomorphically by $\exp_G$ onto an 
open subset of $G$. All Banach--Lie groups 
are locally exponential (\cite[Prop.~IV.1.2]{Ne06}). 
Not every infinite dimensional Lie group has an exponential 
function (\cite[Ex.~II.5.5]{Ne06}), but exponential functions 
are unique whenever they exist. 

If $\pi \: G \to \GL(V)$ is a representation of $G$ on a locally convex 
space $V$, the exponential function 
permits us to associate to each element $x$ of the Lie algebra a 
one-parameter group $\pi_x(t) := \pi(\exp_G tx)$. 
We therefore {\bf assume} in the following that 
$G$ has an exponential function.

\section{Basic facts and definitions} \label{sec:3}

In this section we introduce some basic notation and derive 
some general results for representations of Lie groups on 
locally convex spaces. 

\begin{defn} \label{def:domains}
Let $(\pi, V)$ be a representation of the Lie group $G$ 
(with a smooth exponential function) on the locally convex space $V$. 

(a) We say that $\pi$ is {\it continuous} if the action of $G$ on $V$ 
defined by $(g,v) \mapsto \pi(g)v$ is continuous. 

(b) An element $v \in V$ is a $C^k$-vector, $k \in \N_0 \cup \{\infty\}$,  
if the orbit map $\pi^v \: G \to V, g \mapsto \pi(g)v$ is a 
$C^k$-map. We write $V^k := V^k(\pi)$ for the linear subspace of 
$C^k$-vectors and we say that the representation 
$\pi$ is {\it smooth} if the space $V^\infty$ of smooth vectors 
is dense.  

If $G$ is a Banach--Lie group and $V$ a Banach space, then we write 
$FV^k \subeq V^k$ for the subspace of those $C^k$-vectors whose orbit 
map is also $C^k$ in the Fr\'echet sense. 

(c) For each $x \in \g$, we write 
$$\cD_x := \Big\{ v \in V \: \derat0 \pi(\exp_G tx)v \ \mbox{ exists } \Big\}
$$ 
for the domain of the infinitesimal generator 
$$\oline{\dd\pi}(x)v := \derat0 \pi(\exp_G(tx))v.$$ 
of the one-parameter group $\pi(\exp_G(tx))$. 
We write $\cD_\g := \bigcap_{x \in \g} \cD_x$ 
and $\omega_v(x) := \oline{\dd\pi}(x)v$ for $v \in \cD_\g$. 
Each $\cD_x$ and therefore also $\cD_\g$ are linear subspaces of~$V$, 
but at this point we do not know whether $\omega_v$ is linear (cf.\ 
Theorem~\ref{thm:additive} for a positive answer for Banach--Lie groups). 

(d) We define inductively $\cD_\g^1 := \cD_\g$ and 
$$ \cD_\g^n := 
\{ v \in \cD_\g \: (\forall x \in \g)\, \oline{\dd\pi}(x)v \in 
\cD_\g^{n-1}\} \quad \mbox{ for } \quad n > 1, $$
so that 
$$ \omega_v^n(x_1, \ldots, x_n) := 
\oline{\dd\pi}(x_1)\cdots \oline{\dd\pi}(x_n)v$$ 
is defined for $v \in \cD_\g^n$ and $x_1,\ldots, x_n \in \g$. 
We further put $\cD_\g^\infty := \bigcap_{n \in \N} \cD_\g^n$. 

(e) Let $A$ be an operator with domain $\cD(A)$ on the 
locally convex space $V$. An element in the space 
$\cD^\infty(A) := \bigcap_{n = 1}^\infty \cD(A^n)$ is called a 
{\it smooth vector for $A$}. 

(f) If $V$ is a Banach space, then we say that $(\pi,V)$ 
is {\it locally bounded} if there exists a $\1$-neighborhood 
$U \subeq G$ for which $\pi(U)$ is a bounded set of linear operators 
on $V$. It is easy to see that this requirement implies boundedness of 
$\pi$ on a neighborhood of any compact subset of $G$.
\end{defn}

\begin{rem} \label{rem:cn-vec} 
(a) For every representation $(\pi,V)$ we have 
$$ V^1(\pi) \subeq \cD_\g \quad \mbox{ and } \quad 
V^k(\pi) \subeq \cD_\g^k \quad \mbox{ for } \quad k \in \N. $$
Note that $\omega_v^k$ is continuous and $k$-linear for 
every $v \in V^k(\pi)$. 

(b) By definition, we have $\cD_\g^2 \subeq \cD_\g^1$, so that we obtain 
by induction that $\cD_\g^{n+1} \subeq \cD_\g^n$ for every 
$n \in \N$. 

(c) If $v \in \cD_x$, then $\gamma(t) := \pi(\exp_G(tx))v$ 
is a $C^1$-curve in $V$ with 
$$ \gamma'(t) 
= \pi(\exp_G(tx))\oline{\dd\pi}(x)v, $$ 
so that 
\begin{equation}\
  \label{eq:diff-rel0}
\pi(\exp_G x)v-v = \int_0^1 \pi(\exp_G tx)\oline{\dd\pi}(x)v\, dt. 
\end{equation}
\end{rem}

\begin{lem} \label{lem:a.11c} Suppose that $(\pi, V)$ 
is a continuous representation of the Lie group $G$ 
on $V$. Then a vector $v \in \cD_\g$ 
is a $C^1$-vector if and only if the following two 
conditions are satisfied: 
\begin{description}
\item[\rm(i)] For every smooth curve $\gamma \: [-\eps,\eps]  \to G$ 
with $\gamma(0) = \1$ and $\gamma'(0) = x$, the derivative  
$\frac{d}{dt}|_{t = 0} \pi(\gamma(t))v$ exists and 
equals $\oline{\dd\pi}(x)v$. 
\item[\rm(ii)] $\omega_v \: \g \to V, x \mapsto \oline{\dd\pi}(x)v$ 
is continuous. 
\end{description}

If $G$ is locally exponential, then 
{\rm(i)} follows from {\rm(ii)}, and if 
$G$ is finite dimensional, then {\rm(i)} and {\rm(ii)} 
hold for every $v \in \cD_\g$. 
\end{lem}

\begin{prf} If $\pi^v$ is a $C^1$-map, then 
$\omega_v = T_\1(\pi^v)$ is a continuous linear map satisfying 
$\frac{d}{dt}|_{t = 0}\pi(\gamma(t))v = \omega_v(x)$ for each smooth curve 
in $G$ with $\gamma(0) =\1$ and $\gamma'(0) = x$. 

Suppose, conversely, that (i) and (ii) are satisfied. 
Then the relation 
$$\pi(\gamma(t+h)) = \pi(\gamma(t)) \pi(\gamma(t)^{-1}\gamma(t+h)) $$
implies that for each smooth curve $\gamma$ in $G$, we have 
$$ \frac{d}{dt} \pi(\gamma(t))v = \pi(\gamma(t))\omega_v(\delta(\gamma)_t), $$
where $\delta(\gamma)_t = \gamma(t)^{-1}.\gamma'(t)$ is the left logarithmic 
derivative of $\gamma$ (cf.~\eqref{eq:tang-triv} in Section~\ref{sec:2}). 

We conclude that the orbit map $\pi^v$ 
has directional derivatives in each point $g \in G$, and that its 
tangent map is given by 
$$ T_g(\pi^v)(g.x) = \pi(g)\omega_v(x) $$ 
(cf.\ \eqref{eq:tang-triv}). 
Since $\omega_v$ and the action of $G$ on $V$ are continuous, 
the map $T(\pi^v) \: TG \cong G \times \g \to V$ is continuous, 
i.e., $\pi^v \in C^1(G,V)$. From \cite[Lemma~2.2.14]{GN10} 
it now follows that $\omega_v$ is linear. 

Now we assume, in addition, that $G$ is locally exponential 
and claim that (ii) implies (i). 
Any smooth curve $\gamma$ with $\gamma(0) = \1$ and 
$\gamma'(0) = x$ can be written 
for sufficiently small values of $t \in \R$ as 
$\exp_G(\eta(t))$ with a smooth curve $\eta \: [-\eps,\eps] \to \g$ 
satisfying 
$\eta(0) = 0$ and $\eta'(0) = x$. Therefore it suffices to show that 
for any such curve $\eta$, we have 
$$ \deras0 \pi(\exp_G \eta(s))v = \omega_v(x). $$
First we note that $\eta(0) = 0$ implies that 
$$ \frac{\eta(s)}{s} 
= \frac{1}{s} \int_0^1 \eta'(ts)s\, dt 
= \int_0^1 \eta'(st)\, dt $$
extends by the value $x = \eta'(0)$ to a smooth curve on $[-\eps,\eps]$. 

Next, we derive for $v \in \cD_\g$ from \eqref{eq:diff-rel0} the 
relation 
\begin{equation}
  \label{eq:2}
\pi(\exp_G \eta(s))v-v 
= \int_0^1 \pi(\exp_G t\eta(s))\omega_v(\eta(s))\, dt, 
\end{equation}
which leads to 
\begin{align*}
 \deras0 \pi(\exp_G \eta(s))v
&= \lim_{s \to 0} \int_0^1 \pi(\exp_G t\eta(s))\omega_v(\eta(s)/s)\, dt \\
&= \int_0^1 \omega_v(x)\, dt = \omega_v(x) 
\end{align*}
because, in view of (ii), the integrand in \eqref{eq:2} 
is a continuous function of $(s,t)$, 
so that we may exchange integration and the limit. 
This means that (i) is satisfied. 

If $G$ is finite dimensional, then Goodman's argument in 
\cite[p.~221]{Go70} implies that every $v \in \cD_\g$ 
is a $C^1$-vector, so that (i) and (ii) are satisfied. 
\end{prf}

\begin{lem} \label{lem:a.12} 
If $G$ is locally exponential, then a vector $v \in V$ is a $C^k$-vector 
if and only if $v \in \cD_\g^k$ and the 
maps $\omega_v^n$, $n \leq k$, are continuous and $n$-linear. 
In particular, $v$ is a smooth vector if and only if 
$v \in \cD_\g^\infty$ and all the maps $\omega_v^n$ 
are continuous and $n$-linear. 
\end{lem}

\begin{prf} If $v$ is a $C^k$-vector, then 
$$\omega_v^n(x_1,\ldots, x_k) 
= \frac{\partial^n}{\partial t_1 \cdots\partial t_n}\Big\vert_{t_1 = \cdots = t_n=0}  
\pi(\exp_G(t_1 x_1)\cdots \exp_G(t_n x_n))v $$
is a continuous $n$-linear map for any $n \leq k$. 

Suppose, conversely, that this is the case for some $k \geq 1$. 
Lemma~\ref{lem:a.11c} takes care of the case $k = 1$. 
We may therefore assume that $k > 1$ and that the 
assertion holds for $k - 1$. 
From Lemma~\ref{lem:a.11c} we derive that $v$ is a $C^1$-vector 
with $T(\pi^v)(g.x) = \pi(g) \omega_v^1(x)$. 
We have to show that $T(\pi^v) \: 
TG \cong G \times \g \to V$ 
is a $C^{k-1}$-map. For each fixed $x \in \g$, 
the element $\omega_v^1(x) = \oline{\dd\pi}(x)v$ 
is contained in $\cD_\g^{k-1}$, so that, in view of 
$$ \omega^n_{\omega_v^1(x)}(x_1,\ldots, x_n) 
= \omega_v^{n+1}(x_1, \ldots, x_n,x) \quad \mbox{ for } \quad 
n \leq k -1, $$
 our induction hypothesis implies that 
$\omega_v^1(x)$ is a $C^{k-1}$-vector, hence that $T(\pi^v)$ has 
directional derivatives of all orders $\leq k-1$, and that they are sums 
of terms of the form 
$$ \pi(g)\oline{\dd\pi}(x_1)\cdots \oline{\dd\pi}(x_j)v 
= \pi(g)\omega_v^j(x_1,\ldots, x_j), $$
which are continuous functions 
on $G \times \g^j$, $j \leq k-1$. This proves that 
$T(\pi^v)$ is a $C^{k-1}$-map, and hence that 
$v$ is a $C^k$-vector. 
\end{prf}

\section{A topology on the space of smooth vectors} \label{sec:4}

Throughout this section, we assume that $G$ is a Banach--Lie group 
and that $\|\cdot\|$ is a compatible norm on $\g = \L(G)$ satisfying 
$$ \|[x,y]\| \leq \|x\| \cdot \|y\|, \quad x,y \in \g. $$
We shall define a topology on $V^\infty$ for which 
the action 
$$\sigma \: G \times V^\infty \to V^\infty, \quad 
(g,v) \mapsto \pi(g)v $$ 
is smooth. 

Let 
$$\dd \pi \: \g \to \End(V^\infty), \quad 
\dd\pi(x)v := \derat0 \pi(\exp tx)v $$
denote the derived  action of $\g$ on $V^\infty$. 
That this is indeed a representation of $\g$ follows by observing that 
the map $V^\infty \to C^\infty(G,V), v \mapsto \pi^v$ 
intertwines the action of $G$ with the right translation action on 
$C^\infty(G,V)$, and this implies that the derived action 
corresponds to the action of $\g$ on $C^\infty(G,V)$ by left invariant
 vector fields (cf.\ \cite[Rem.~IV.2]{Ne01a} for details). 

\begin{defn} \label{def:4.1} (a) For $n \in \N_0$, let 
$\Mult^n(\g,V)$ be the space of continuous $n$-linear maps 
$\g^n \to V$ (for $n = 0$ we interpret this as the constant maps) 
and write $\cP(V)$ for the set of continuous seminorms on $V$. 
The space $\Mult^n(\g,V)$ carries a natural locally convex topology defined by 
the seminorms 
$$ p(\omega) := \sup \{ p(\omega(x_1,\ldots, x_n)) \: 
\|x_1\|, \ldots, \|x_n\| \leq 1\}, \quad 
p \in \cP(V). $$
Note  that $p(\omega)$ is the smallest constant $c \geq 0$ for which 
we have the estimate 
\begin{equation}\label{eq:derivest0}
p(\omega(x_1, \ldots, x_n)) \leq c \|x_1\| \cdots \|x_n\|\quad \mbox{ for } 
\quad x_1, \ldots, x_n \in \g. 
\end{equation}

(b) To topologize $V^\infty$, we define for each $n \in \N_0$ a map 
$$ \Psi_n \: V^\infty \to \Mult^n(\g,V), \quad 
\Psi_n(v)(x_1,\ldots, x_n) := \dd\pi(x_1) \cdots \dd\pi(x_n)v. $$
That $\Psi_n(v)$ defines a continuous $n$-linear map 
follows from the smoothness of the orbit map 
$\pi^v \: G \to V, g \mapsto \pi(g)v$ and the fact that 
$\Psi_n(v)$ is obtained by $n$-fold partial derivatives in $\1$. 
We thus obtain an injective linear map 
$$ \Psi \: V^\infty \to \prod_{n \in \N_0} \Mult^n(\g,V), \quad 
v \mapsto (\Psi_n(v))_{n \in \N_0}, $$
and define the topology on $V^\infty$ such that 
$\Psi$ is a topological embedding. 
This means that the topology on $V^\infty$ is defined by 
the seminorms 
$$ p_n(v) := \sup \{ p(\dd\pi(x_1) \cdots \dd\pi(x_n)v) \: 
x_i \in \g, \|x_i\|\leq  1\}, \quad p \in \cP(V), n \in \N_0.$$ 
Note  that $p_n(v) = p(\Psi_n(v))$ 
is the smallest constant $c \geq 0$ for which 
we have the estimate 
\begin{equation}\label{eq:derivest}
p(\dd\pi(x_1) \cdots \dd\pi(x_n)v)
 \leq c \|x_1\| \cdots \|x_n\|, \quad 
x_1, \ldots, x_n \in \g. 
\end{equation}
We endow $V^\infty$ with the locally convex topology defined by 
the seminorms $p_n$, $n \in \N_0$, $p \in \cP(V)$.  
\end{defn}

\begin{lem} \label{lem:4.1.1} 
 The bilinear  map $\dd \pi \: \g \times V^\infty \to  V^\infty, (x,v) 
\mapsto \dd\pi(x)v $ is continuous.  
\end{lem}

\begin{prf} This follows from 
$p_n(\dd\pi(x)v) \leq p_{n+1}(v)\|x\|$, which is a consequence 
of \eqref{eq:derivest}. 
\end{prf}

\begin{lem} \label{lem:4.1.2} 
 The group $G$ acts by continuous linear operators on 
$V^\infty$. More precisely, we have 
\begin{equation}
  \label{eq:esti2}
 p_n(\pi(g)v) \leq (p \circ \pi(g))_n(v) \|\Ad(g)^{-1}\|^n. 
\end{equation}
\end{lem}

\begin{prf} For $x_1, \ldots, x_n \in \g$, we have 
$$ \dd\pi(x_1) \cdots \dd\pi(x_n) \pi(g) v
=  \pi(g) \dd\pi(\Ad(g^{-1})x_1) \cdots \dd\pi(\Ad(g^{-1})x_n) v, $$
and, in view of \eqref{eq:derivest}, this implies \eqref{eq:esti2}.   
\end{prf}

\begin{thm} \label{thm:4.1.3} 
If $(\pi, V)$ is a representation of the Banach--Lie group 
$G$ on the locally convex space $V$ defining a continuous action 
of $G$ on $V$, then the 
action $\sigma(g,v) := \pi(g)v$ of $G$ on $V^\infty$ is smooth. 
\end{thm}

\begin{prf} First we show that $\sigma$ is continuous. 
In view of Lemma~\ref{lem:4.1.2}, it  
suffices to show that $\sigma$ is continuous in 
$(\1,v)$ for each $v \in V^\infty$. 
For $w \in V^\infty$, the relation 
$$ \pi(g)w - v = \pi(g)(w-v) + \pi(g)v - v $$ 
permits us to break the argument into a proof for the 
continuity of $\sigma$ in $(\1,0)$ and the continuity of the 
orbit map $\pi^v \: G \to V^\infty$. 

First we use \eqref{eq:esti2} to obtain
$$ p_n(\pi(g)w) 
\leq (p  \circ \pi(g))_n(w)\|\Ad(g)^{-1}\|^n. $$
Since $U^p := \{ w \in V \: p(w) \leq 1\}$ is a $0$-neighborhood 
in $V$, the continuity of the action of $G$ on $V$ implies the 
existence of a $\1$-neighborhood $U_G \subeq G$ and a 
continuous seminorm $q$ on $V$ with $\pi(U_G) U^q \subeq U^p$. 
This means that $p(\pi(g)v) \leq q(v)$ for $v \in V$ and $g \in U_G$, 
i.e., $p \circ \pi(g) \leq q$. We  thus obtain  
for $g \in U_G$ the estimate 
$$ p_n(\pi(g)(w)) 
\leq q_n(w)\|\Ad(g)^{-1}\|^n,  $$
and since $\Ad \: G \to \GL(\g)$ is locally bounded, 
$\sigma$ is continuous in $(\1,0)$. 

To verify the continuity of the orbit maps 
$\pi^v$, $v \in V^\infty$, we consider the smooth function 
$f \: G \times G \to V, (g,h) \mapsto \pi(g)\pi(h)v$ and 
observe that $f^h(g) := f(g,h) = \pi^{\pi(h)v}(g)$. 
Using the family $\cP(V)$ of all continuous seminorms on $V$, 
we embed $V$ into the topological product 
$\prod_{p \in \cP(V)} V_p$, where $V_p$ is the Banach space 
obtained by completing $V/p^{-1}(0)$ with respect to the 
norm on this space induced by $p$. We write $v \mapsto [v]_p$ for the 
corresponding quotient map. Then the smoothness 
of the map $f$ is equivalent to the smoothness of the 
component mappings $f_p \: G \times G \to V_p, (g,h) \mapsto [f(g,h)]_p$, and 
since smoothness for Banach space-valued maps implies smooth dependence 
of higher derivatives (cf.\ \cite[Thm.~I.7]{Ne01a}, \cite{KM97}, \cite{GN10}), 
it follows that the 
Banach space-valued maps 
$$ f_p^n \: G \to \Mult^n(\g,V_p), \quad h \mapsto \Psi_n(\pi(h)v)_p $$
are smooth. Here we use that for each open subset $U$ in a 
Banach space $E$ and a smooth map $F \: U \to V_p$, the map 
$F^n \: U \to \Mult^n(E,V_p)$, defined by 
$$ F^n(x)(v_1,\ldots, v_n) := 
(\partial_{v_1} \cdots \partial_{v_n} F)(x), $$ 
is smooth.   This implies that the corresponding map 
$$ f_p \: G \to \prod_{n \in \N_0} \Mult^n(\g,V_p), 
\quad h \mapsto \Psi(\pi(h)v)_p = \Psi(\pi^v(h))_p$$ 
is smooth, and this in turn means that the orbit map 
$\pi^v\: G \to V^\infty$ is smooth, hence in particular continuous. 
This completes the proof of the continuity of $\sigma$. 

The preceding argument already implies that the partial derivatives 
of the action map $\sigma$ exist and that they are given by 
$$ T_{(g,v)}(\sigma)(g.x,w) = \pi(g)\dd\pi(x)v + \pi(g)w. $$
From the continuity of $\sigma$ and the action of $\g$ on $V^\infty$ 
(Lemma~\ref{lem:4.1.1}), it now follows that $T\sigma$ is continuous, so that 
$\sigma$ is actually $C^1$. Iterating this argument, we see that whenever 
$\sigma$ is $C^n$, then $T\sigma$ also is $C^n$, and this 
shows that $\sigma$ is smooth. 
\end{prf}

\begin{cor}
  \label{cor:hilsmooth} 
If $G$ is a Banach--Lie group and 
$(\pi, \cH)$ a continuous unitary representation of $G$ on $\cH$, 
then the induced action of $G$ on $\cH^\infty$ is smooth. 
\end{cor}

It is instructive to compare our topology on 
$V^\infty$ with the standard construction for finite dimensional 
Lie groups: 

\begin{prop} \label{prop:findim} Suppose that $G$ is finite dimensional, 
$V$ is a locally convex space, and 
$C^\infty(G,V)$ is endowed with the smooth compact open topology, i.e., 
the topology of uniform convergence of all derivatives on compact 
subsets. Then for any continuous representation $(\pi,V)$ of $G$, 
the injection 
$$ \eta \: V^\infty \to C^\infty(G,V), \quad v \mapsto \pi^v $$
is a topological embedding 
whose range is the subspace of smooth maps $f \: G \to V$ which are 
equivariant in the sense that 
$$ f(gh) = \pi(g) f(h) \quad \mbox{ for } \quad g,h \in G. $$
\end{prop}

\begin{prf} We recall that the smooth compact open topology is 
defined by the property that the 
embedding $C^\infty(G,V) \into \prod_{n \in \N_0} C(T^n G, T^n V), 
f \mapsto (T^n(f))_{n \in \N_0}$, 
defined by the tangent maps, is a topological embedding, where the 
factors on the right are endowed with the compact open topology. 
It is easy to see that this is a locally convex topology for which 
the action of $G$ on 
$C^\infty(G,V)$ by right translations is smooth 
(cf. \cite[Sect.~III]{Ne01a}). Here the key observation is that the 
smoothness of the map 
$$ G \times C^\infty(G,V) \to C^\infty(G,V), \quad 
(g,f) \mapsto f \circ \rho_g $$
follows from the smoothness of the corresponding map 
$$ G \times C^\infty(G,V) \times G \to V, \quad 
(g,f,x) \mapsto f(xg) $$ 
(\cite[Lemma~A.2]{NW08}). This in turn 
follows from the smoothness of the evaluation map 
$\ev \: C^\infty(G,V) \times G \to V$, whose smoothness follows 
from the smoothness of 
$\id_{C^\infty(G,V)}$ (cf.\ \cite[Lemma~A.3(1)]{NW08}; resp., 
\cite[Prop.~12.2]{Gl04}).\begin{footnote}
{In the setting of convenient calculus, one finds similar arguments 
for the smoothness of the $G$-action on $V^\infty$ 
in \cite[Thm.~5.2]{Mi90}.}
\end{footnote}

If $(\pi,V)$ is a continuous representation of 
$G$ on $V$, then $\eta(v) := \pi^v$ injects 
$V^\infty$ into $C^\infty(G,V)$ and the image of $\eta$ 
clearly is the subspace of equivariant map because 
any such map $f$ satisfies $f = \pi^{f(\1)}$. 
Since $\im(\eta)$ is a closed subspace invariant under right 
translation, we obtain by restriction a smooth 
action on this space. 

For $d = \dim \g$, we also have a topological 
isomorphism $\Mult^p(\g,V) \cong V^{pd}$ by evaluating 
$p$-linear maps on $p$-tuples of elements of a fixed basis of~$\g$. 
Therefore the topology on $V^\infty$ 
is defined by the seminorms 
$$ v \mapsto p(\dd\pi(x_1)\cdots \dd\pi(x_n)v),\quad 
p \in \cP(V), $$
where $x_1, \ldots, x_n$ are $n$ elements of a fixed basis 
of $\g$ (Definition~\ref{def:4.1}). 

Since $\g$ acts on $C^\infty(G,V)$ by continuous linear maps, 
all the operators $\dd\pi(x)$ are continuous with respect to the 
smooth compact open topology, and therefore the evaluation map 
$$ \ev_\1 \: \im(\eta) \to V^\infty$$
is a continuous bijection. 
Conversely, the smoothness of the action of $G$ on $V^\infty$ 
(Theorem~\ref{thm:4.1.3}) implies that the linear map $\eta$ 
is continuous (\cite[Lemma~A.2]{NW08}), hence a topological 
embedding. 
\end{prf}

\begin{rem} (a) If $G$ is finite dimensional with countably 
many connected components, which is equivalent to being a countable 
union of compact subsets, and $V$ is a Fr\'echet space, 
then $C^\infty(G,V)$ is also Fr\'echet, and therefore 
$V^\infty$ is a Fr\'echet space. 

(b) For a unitary representation $(\pi,  \cH)$ of a finite dimensional 
Lie group $G$, Goodman uses in \cite{Go69} Sobolev space techniques 
to show that the space $\cH^\infty$ of smooth vectors is 
the intersection of the spaces 
$\cD^\infty(\oline{\dd\pi}(x_j)), j =1,\ldots, n$, where 
$x_1, \ldots, x_n$ is a basis for $\g$ (\cite[Thm.~1.1]{Go69}; 
see also \cite[Thm.~10.1.9]{Sch90}). 
This implies in particular that $\cH^\infty$ is complete 
and that the topology on this space is compatible 
with our construction 
(cf.\ Proposition~\ref{prop:findim}). 
\end{rem}

It seems that the chances that the results in this section 
extend to some classes of Fr\'echet--Lie groups are not very high, 
as the following two examples show. 

\begin{ex}\label{ex:4.7} For the Fr\'echet--Lie group 
$G = (\R^\N,+)$ (endowed with the product topology), 
we consider the continuous unitary representation 
on $\cH = \ell^2(\N,\C)$ given by 
$(\pi(g)x)_n = e^{ig_n} x_n$. Then 
$$ \cD_\g = \{ x \in \cH \: (\forall y \in \R^\N) (x_n y_n) \in \cH \} 
= \Spann \{ e_n \: n \in \N\} = \C^{(\N)} $$ 
is a countably dimensional vector space. 
Since $\cD_\g$ is spanned by eigenvectors, we see that 
$\cD_\g = \cH^\infty$. 

We claim that for no locally convex topology on $\cH^\infty$ the 
bilinear map 
$$ \beta \: \g \times \cH^\infty \to \cH, \quad 
(x,v) \mapsto \dd\pi(x)v $$
is continuous, which implies in particular that the action of 
$G$ on $\cH^\infty$ is not $C^1$, hence in particular not smooth. 
To verify our claim, let $B \subeq \cH$ denote the closed unit ball. 
If $\beta$ is continuous, then there exist $0$-neighborhoods 
$U_\g \subeq \g$ and $U \subeq \cH^\infty$ with 
$\beta(U_\g \times U) \subeq B$. Next we observe that 
$U_\g$ contains a subspace of the form $\R^M$, 
$M := \N \setminus \{ 1,\ldots, n\}$. In particular 
$\R e_{n+1} \subeq U_\g$. Since we also have $\eps e_{n+1} \in U$ for 
some $\eps > 0$, we arrive at the contradiction 
$\R \cdot \eps e_{n+1} = \R e_{n+1} \subeq B$. 

This proves that Theorem~\ref{thm:4.1.3} does not generalize 
to unitary representations of general Fr\'echet--Lie groups. 
\end{ex}

\begin{ex}\label{ex:4.8} We have seen in Proposition~\ref{prop:findim} 
that, for finite dimensional Lie groups, we obtain the natural topology 
on $V^\infty$ by embedding it into the space 
$C^\infty(G,V)$, endowed with the smooth compact open topology. 
In this case the smoothness of the $G$-action on 
$C^\infty(G,V)$ yields the smoothness on the invariant subspace 
$V^\infty$. For any infinite dimensional Lie group $G$, 
the smooth compact open topology still makes sense, but, as the 
following example shows, this topology is too weak to guarantee the 
continuity of the $G$-action. 

To substantiate this claim, we consider a locally convex space 
$G$, considered as a Lie group, and the space 
$E := \Aff(G,\R)$ of affine real-valued functions on $G$.  
Then $G$ acts on $E$ by $(\pi(g)f)(x) := f(x + g)$. 
Identifying $E$ with $\R \times G'$, we see that for 
$f = c + \alpha$ ($c \in \R, \alpha \in G'$) we have 
$$ \pi(g)(c + \alpha) = \alpha(g) + c + \alpha. $$
Therefore $\pi$ is continuous with respect to a locally convex 
topology on $E$, resp., $G'$, if and only if the evaluation map 
$$ G' \times G \to \R, \quad (\alpha,x) \mapsto \alpha(x) $$
is continuous. 

The smooth compact open topology on $C^\infty(G,\R)$ 
corresponds on $G'$ to the topology of uniform convergence 
on compact subsets of $G$. The continuity of the evaluation map 
with respect to this topology is equivalent to the existence 
of a compact subset $C \subeq G$ for which the closed convex 
hull $\oline{\conv(C)}$ is a $0$-neighborhood. 
If $V$ is complete, then the precompactness of 
$\oline{\conv(C)}$ (\cite[Ch.~II, \S 4, no.~2, Prop.~3]{Bou07})  
implies that $G$ has a compact $0$-neighborhood, hence that 
it is finite dimensional. 

A similar argument shows that, for the finer topology of uniform 
convergence on bounded subsets of $G'$, the evaluation map 
is continuous if and only if $G$ has a bounded $0$-neighborhood 
which means that it is a normed space. 
\end{ex}

\subsection*{Distribution vectors} 

An element of the topological dual space $V^{-\infty} := (V^\infty)'$ 
is called a {\it distribution vector}. The main property of these 
functionals is that they correspond to $G$-morphisms 
$V^\infty \to C^\infty(G)$. Note that 
$(g.\alpha) := \alpha \circ \pi(g)^{-1}$ defines a natural action 
of $G$ on $V^{-\infty}$. 

\begin{lem} \label{lem:distvec} We have an injective map 
$$ \Phi \: V^{-\infty} \to \Hom_G(V^\infty, C^\infty(G)), 
\quad \Phi(\alpha)(v)(g) := \alpha(\pi(g)^{-1}v), $$
where $G$ acts on $C^\infty(G)$ by $(g.f)(h) := f(g^{-1}h)$. 
The range of $\Phi$ consists of all those $G$-morphisms 
$\phi \: V^\infty \to C^\infty(G)$ for which the composition 
$\ev_\1 \circ \phi$ is continuous, i.e., an element of 
$V^{-\infty}$. 
\end{lem}

\begin{prf} The smoothness of the functions $\Phi(\alpha)(v)$ 
follows from the smoothness of the $G$-action on $V^\infty$, and 
the equivariance of $\Phi(\alpha)$ is immediate from the definition. 

For $\phi := \Phi(\alpha)$ we have $\alpha = \ev_\1 \circ \phi 
\in V^{-\infty}$. If, conversely, 
$\phi \in \Hom_G(V^\infty, C^\infty(G))$ is such that 
$\alpha := \ev_\1 \circ \phi$ is continuous, we have 
$$ \phi(v)(g) 
= (g^{-1}.\phi(v))(\1) 
= \phi(\pi(g)^{-1}v)(\1) 
= (\ev_\1 \circ \phi)(\pi(g)^{-1}v)
= \Phi(\alpha)(v)(g), $$
i.e., $\phi = \Phi(\alpha)$. 
\end{prf}

\begin{rem} (a) 
For finite dimensional Lie groups, compactly supported  
distributions are defined as 
the elements of the topological dual 
$C^\infty(G)'$ of the space of smooth functions. 
Since the $G$-action on $C^\infty(G)$ is smooth 
(cf.\ the proof of Proposition~\ref{prop:findim}), compactly supported 
distributions are distribution vectors in the sense defined above. 

(b) We do not know if there exists a locally convex topology 
on $C^\infty(G)$ for which the left or right translation action 
is smooth and therefore it is not clear what a good concept 
of a distribution on a Banach--Lie group is. 
However, the topology on $V^\infty$ for a continuous 
representation $(\pi,V)$ provides a natural concept of a 
distribution vector. 
\end{rem}

In harmonic analysis on homogeneous spaces $G/H$, distribution 
vectors of unitary representations play an important role. 
If $H \subeq G$ is a Lie subgroup, so that the quotient 
$G/H$ carries a manifold structure for which the 
quotient map $q \: G \to G/H, g \mapsto gH,$ is a submersion, then we may 
identify $C^\infty(G/H)$ with the subspace of $C^\infty(G)$ 
consisting of all functions constant on the cosets $gH$. 
For the map $\Phi$ in Lemma~\ref{lem:distvec} we immediately see that 
$$ \Phi(\alpha)(V^\infty) \subeq C^\infty(G/H) $$
is equivalent to the invariance of $\alpha$ under $H$. 
Therefore the space $(V^\infty)^{H}$ of $H$-invariant 
distribution vectors parametrizes the 
$G$-morphisms $\phi \: V^\infty \to C^\infty(G/H)$ 
which are continuous in the very weak sense that their 
composition with point evaluations is continuous (Lemma~\ref{lem:distvec}, 
see also \cite[p.~137]{HS94}). 
This situation is of particular interest if these embeddings 
are essentially unique: 
A pair $(G,H)$ of a Banach--Lie group $G$ and a subgroup $H$ 
is called a {\it generalized Gelfand pair} if 
$$ \dim(\cH^{-\infty})^H \leq 1 $$ 
holds for every irreducible continuous unitary representation $(\pi,\cH)$ of 
$G$. For finite dimensional Lie groups there exists a 
rich theory of generalized Gelfand pairs with many 
applications in harmonic analysis (cf.\ \cite{vD09}). 
It is a very worthwhile project to explore this concept 
also for Banach--Lie groups, e.g., for automorphism groups  
of infinite dimensional semi-Riemannian symmetric spaces 
such as hyperboloids 
$$ X = \{ (t,v) \in \R \times \cH \: \|v\|^2 - t^2 = 1 \}, $$
where $\cH$ is a real Hilbert space 
(cf.\ \cite{Fa79} for the case $\dim \cH < \infty$).

\section{Smooth vectors in Banach spaces} \label{sec:5}

In this section we assume that $G$ is a Banach--Lie group 
and that $V$ is a Banach space. Our goal 
is to show that the space $V^\infty$ is complete, hence a Fr\'echet 
space.

The following proposition generalizes the well-known fact that 
separately continuous bilinear maps on Fr\'echet spaces are continuous. 

\begin{prop} {\rm(Continuity criterion)} \label{prop:concrit} 
Let $X$ be a first countable topological space, 
$F$ a Fr\'echet space and $V$ a topological vector space. 
Let $\alpha \: X \to \cL(F,V)$ (the space of continuous linear operators 
$F \to V$) 
be a map such that, for each $f \in F$, the map 
$$ \alpha_f \: X \to V, \quad x \mapsto \alpha(x)f $$
is continuous. Then the map 
$\hat\alpha \: X \times F \to V, (x,f) \mapsto \alpha(x)f$
is continuous. 
\end{prop}

\begin{prf} Since $X$ is first countable, the same holds 
for the product space $X \times F$, so that we only have
to show that, for any sequence $(x_n, f_n) \to (x_0, f_0)$ in $X \times F$, 
we have $\alpha(x_n)f_n \to \alpha(x_0)f_0$. 
Our assumption implies that 
the sequence $\alpha(x_n)$ of linear maps 
converges pointwise to $\alpha(x_0)$. In particular, 
for each $f \in F$, the sequence $\alpha(x_n)f$ in $V$ is bounded. 
Now the Banach--Steinhaus 
Theorem (\cite[Thm.~2.6]{Ru73}) implies that the sequence 
$\alpha(x_n)$ is equicontinuous, which leads to 
$\alpha(x_n)(f_n - f_0) \to 0.$
We thus obtain 
$$ \alpha(x_n)f_n - \alpha(x_0)f_0 
=  \alpha(x_n)f_n - \alpha(x_n)f_0 
 + \alpha(x_n) f_0  - \alpha(x_0)f_0 \to 0. $$
This proves the continuity of $\hat\alpha$. 
\end{prf}

\begin{lem} \label{lem:4.1} 
Let $V$ be a Banach space, $G$ be a topological 
group and $\pi \: G \to \GL(V)$ be a homomorphism. 
Then the following are equivalent:
\begin{description}
\item[\rm(i)] The linear action $\sigma  \: G \times V \to V, (g,v) \mapsto \pi(g)v$ is continuous. 
\item[\rm(ii)] $\sigma$ is continuous in $(\1,v)$ for every $v \in V$. 
\item[\rm(iii)] $\pi$ is locally bounded and $\pi$ is 
strongly continuous, i.e., all orbit maps 
$\pi^v(g) = \sigma(g,v)$ are continuous. 
\end{description}
If, in addition, $G$ is metrizable, then {\rm(i)}-{\rm(iii)} are 
equivalent to 
\begin{description}
\item[\rm(iv)] $\pi$ is strongly continuous. 
\end{description}
\end{lem}

\begin{prf} (i) $\Rarrow$ (ii) is trivial. 

(ii) $\Rarrow$ (iii): If $\sigma$ is continuous in 
$(\1,0)$,   then there exists a $\1$-neighborhood $U\subeq G$ 
and a ball $B_\eps(0)$ in $V$ with 
$\pi(U) B_\eps(0) \subeq B_1(0)$. This implies 
that $\|\pi(g)\| \leq \eps^{-1}$ for every $g \in U$. 
The continuity in the pairs $(\1,v)$ implies the continuity of the 
orbit maps. 

(iii) $\Rarrow$ (i) follows from 
$\|\pi(g)v - \pi(h)w\| 
\leq \|\pi(g)\| \|v - w \| + \|\pi(g)w - \pi(h)w \|.$

(iv) $\Rarrow$ (i): If, in addition, $G$ is metrizable, then it is first 
countable, so that Proposition~\ref{prop:concrit} shows that 
(iv) implies (i). 
\end{prf} 

\begin{rem} \mlabel{rem:5.3} If $\pi \: G \to \GL(V)$ is a representation 
of the topological group $G$ on the Banach space $V$ which 
is locally bounded, then the subspace 
$V^0$ of {\it continuous vectors}, i.e., of those 
$v \in V$ for which the orbit map $\pi^v \: G \to V, g \mapsto \pi(g)v$, 
is continuous is closed. In fact, suppose that $v_n \to v$ holds 
for some sequence $v_n \in V^0$ and that 
$U$ is a neighborhood of $g_0 \in G$ on which $\pi$ is bounded. 
Then $\pi^{v_n} \to \pi^v$ holds uniformly on $U$. Therefore 
the continuity of the maps $\pi^{v_n}$ implies the 
continuity of $\pi^v$ on $U$. Since $g_0$ was arbitrary, it follows that 
$v \in V^0$, and hence that $V^0$ is closed. Since $V^0$ is 
$G$-invariant, we obtain a locally bounded representation 
$\pi^0 \: G \to \GL(V^0)$, and its continuity follows from 
Lemma~\ref{lem:4.1}. 
\end{rem}

In view of the preceding lemma, the continuity of 
the action of $G$ on $V$ implies that $\pi$ is locally bounded. 
The following proposition can also be derived from the results 
in Section~\ref{sec:10} (cf.\ Remark~\ref{rem:10.6}). 

\begin{prop} \label{prop:5.2} 
If $(\pi,V)$ is a continuous representation of the Banach--Lie group $G$ 
on the Banach space $V$, then $V^\infty$ is complete, 
i.e., a Fr\'echet space. 
\end{prop}

\begin{prf} First we note that $\prod_{n \in \N_0} \Mult^n(\g,V)$ 
is a countable product of Banach spaces, hence a Fr\'echet space. 
Therefore the completion $\hat V^\infty$ of the topological 
vector space $V^\infty$ can be identified with the closure 
of the subspace $\Psi(V^\infty)$. Since $V$ is complete, we have a 
continuous linear map 
$\iota \: \hat V^\infty \to V$, extending the inclusion 
$V^\infty \into V$, and, with respect to the realization 
of $\hat V^\infty$ as the closure of the image of $\Psi$, 
this map is given by $\iota((\alpha_n)_{n \in \N_0}) = \alpha_0$. 
In the course of the proof we shall see that $\iota$ is injective 
and that its range coincides with $V^\infty$. This implies that 
$\hat V^\infty$ is not larger than $V^\infty$, so that 
$V^\infty$ is complete. 

Let $(v_n)_{n \in\N}$ be a sequence in $V^\infty$ for which 
$\Psi_k(v_n)$ converges to some $\omega^k \in \Mult^k(\g,V)$ for each 
$k \in \N_0$. We have to show that $v := \omega^0$ is a smooth vector 
and that $\omega^k = \Psi_k(v)$ for each $k$. 
Since the continuous linear maps 
$\dd\pi(x) \: V^\infty \to V^\infty$ extend to continuous linear maps 
$$\hat\dd\pi(x) \: \hat V^\infty \to \hat V^\infty$$ 
on the completion, the convergence $v_n \to v$ in 
$\hat V^\infty$ implies that 
$$ \dd\pi(x_1)\cdots \dd\pi(x_k) v_n 
\to \hat\dd\pi(x_1)\cdots \hat\dd\pi(x_k) v, $$
and hence that 
\begin{equation}
  \label{eq:omegarel}
\omega^k(x_1,\ldots, x_k) =  \hat\dd\pi(x_1)\cdots \hat\dd\pi(x_k) v 
\in \hat V^\infty.
\end{equation}
This proves in particular that $v = 0$ implies 
$\omega^k = 0$ for each $k > 0$, and hence that the map 
$\iota \: \hat V^\infty \to V$ is injective. 

For $x \in \g$, we obtain with \eqref{eq:diff-rel0} 
that 
$$ \pi(\exp_G x)v_n -v_n 
= \int_0^1 \pi(\exp_G tx)\dd\pi(x) v_n\, dt $$
for each $n$. Since  $\dd\pi(x)v_n \to \omega^1(x)$ 
and the linear map 
$w \mapsto \int_0^1 \pi(\exp_G(tx))w\, dt$ 
is continuous ($\pi$ is locally bounded), we obtain 
\begin{equation} 
  \label{eq:diff-rel}
\pi(\exp_G x)v-v = \int_0^1 \pi(\exp_G tx)\omega^1(x)\, dt, 
\end{equation}
which leads to 
$$ \deras0 \pi(\exp_G sx)v 
= \lim_{s \to 0} \int_0^1 \pi(\exp_G stx)\omega^1(x)\, dt
= \int_0^1 \omega^1(x)\, dt = \omega^1(x). $$
This implies that $v \in \cD_\g$ with $\omega_v = \omega^1$. 
In particular, $\omega_v$ is continuous and linear. 
Since Banach--Lie groups are locally exponential, 
Lemma~\ref{lem:a.11c} 
now implies that $v$ is a $C^1$-vector. 

To see that $\pi^v$ is $C^2$, we have to show that 
$T(\pi^v)$ is $C^1$. From \eqref{eq:omegarel} we recall 
that, for each $x \in \g$, 
$\omega^1(x) = \hat\dd\pi(x)v \in \hat V^\infty$ is a $C^1$-vector 
by the preceding argument. 
Therefore 
$$T(\pi^v) \: TG \cong G \times \g \to V, \quad (g,x)
 \mapsto \pi(g)\hat\dd\pi(x)v $$ 
has directional derivatives given by 
$$ T_{g,x}T(\pi^v)(g.y,w) = \pi(g)\hat\dd\pi(y) \hat \dd\pi(x)v 
+ \pi(g)\hat\dd\pi(y)v. $$
We have already seen above that the second term is continuous,  
and the continuity of the first term follows from the continuity of the 
action of $G$ on $V$ and the continuity of the bilinear map 
$$ \omega^2(y,x) = \hat\dd\pi(y) \hat \dd\pi(x)v. $$
This proves that each $\pi^v$ is $C^2$. 
Iterating this argument implies that $\pi^v$ has 
directional derivatives of any order $k$, and that they are sums 
of terms of the form 
$$ \pi(g)\hat\dd\pi(x_1)\cdots \hat\dd\pi(x_j) 
= \pi(g)\omega^j(x_1,\ldots, x_j), $$
which are continuous on $G \times \g^j$. Therefore 
$v$ is a $C^k$-vector for any $k$, hence smooth, 
and thus $\hat V^\infty \subeq V^\infty$ implies that 
$V^\infty$ is complete. 
\end{prf}

\section{An applications to $C^*$-dynamical systems}  \label{sec:6}

In this section we show that in the special situation 
where a Banach--Lie group $G$ acts by automorphisms 
on a unital $C^*$-algebra $\cA$, i.e., for $C^*$-dynamical 
systems, the Fr\'echet space $\cA^\infty$ is a continuous 
inverse algebra, i.e., its unit group is open and the 
inversion is a continuous map. 

\begin{defn} \label{def:a.7} 
(a) A locally convex unital algebra $\cA$ is called a {\it continuous 
inverse algebra} if its group of units $\cA^\times$ is open 
and the inversion map $\eta \: \cA^\times \to \cA, a \mapsto a^{-1}$ 
is continuous. 

(b) Let $G$ be a topological group and $\cA$ be a $C^*$-algebra. 
A {\it $C^*$-dynamical system} is a triple 
$(\cA, G, \alpha)$, where $\alpha \: G \to \Aut(\cA), g \mapsto \alpha_g$, 
is a homomorphism defining a continuous action of $G$ on $\cA$. 
\end{defn}

\begin{thm}\label{thm:a.8}
  \label{cor:cstarsmooth} 
If $G$ is a Banach--Lie group and 
$(\cA, G, \alpha)$ a $C^*$-dynamical system with a unital 
$C^*$-algebra $\cA$, then 
the space $\cA^\infty$ of smooth vectors is a subalgebra, 
which is a topological algebra with respect to its natural topology, 
and the action of $G$ on the Fr\'echet space $\cA^\infty$ is smooth. 

If, in addition, $\cA$ is unital, then $\cA^\infty$ is a continuous 
inverse algebra. 
\end{thm}

\begin{prf} In view of Theorem~\ref{thm:4.1.3}, the action of 
$G$ on $\cA^\infty$ is smooth, and it is a Fr\'echet space 
by Proposition~\ref{prop:5.2}. It therefore remains to 
verify that $\cA^\infty$ is a topological algebra, i.e., 
a subalgebra on which the product and the involution are continuous. 

If $a, b \in \cA^\infty$, then the orbit map of $ab$ is given by 
$\alpha_g(ab) = \alpha_g(a) \alpha_g(b)$, and since the multiplication 
in $\cA$ is continuous bilinear, hence smooth, $ab \in \cA^\infty$. 
We also derive from the smoothness of the inversion in $\cA^\times$ 
that if $a \in \cA^\times \cap \cA^\infty$, then also 
$a^{-1} \in \cA^\infty$. 

The continuity of the multiplication in $\cA^\infty$ follows 
from the fact that the operators $\dd\pi(x) \: \cA^\infty \to \cA^\infty$ 
are derivations. If we write for 
a subset 
$$S =  \{ i_1 < i_2 < \ldots < i_k\} \subeq \{1,\ldots, n\}$$ 
$$ \dd\pi(x_S) := \dd\pi(x_{i_1}) \cdots \dd\pi(x_{i_k}), $$
then 
\begin{equation}
  \label{eq:leibniz}
\dd\pi(x_1) \cdots \dd\pi(x_n)(ab) 
= \sum_S \dd\pi(x_S)(a) \dd\pi(x_{S^c})(b).  
\end{equation}
In view of the submultiplicativity of the norm, this leads to 
$$ \|\dd\pi(x_1) \cdots \dd\pi(x_n)(ab)\| 
\leq \sum_S \|\dd\pi(x_S)(a)\| \cdot \|\dd\pi(x_{S^c})(b)\|,  $$
so that the seminorm $p(a) := \|a\|$ satisfies 
$p_n(ab) \leq \sum_S p_{|S|}(a) p_{|S^c|}(b).$ 
Since the sequence $(p_n)_{n \in \N}$ of seminorms defines 
the topology on $\cA^\infty$, it follows that the multiplication 
in $\cA$ is continuous. 

The continuity of the involution on $\cA^\infty$ follows from 
the fact that the operators $\dd\pi(x)$ commute with $*$: 
$$ \dd\pi(x_1) \cdots \dd\pi(x_n)(a^*) =  \dd\pi(x_1) \cdots 
\dd\pi(x_n)(a)^*,$$
which leads to $p_n(a^*) = p_n(a)$ for $a \in \cA^\infty$. 

Since the inclusion $\cA^\infty \into \cA$ is continuous, 
the unit group 
$$(\cA^\infty)^\times = \cA^\infty \cap \cA^\times$$ 
is open. To prove the continuity of the inversion, 
we show by induction that it is continuous with respect 
to the seminorms $p_n$, $n \leq N$. 
For $N = 0$, the assertion follows from the continuity 
of the inversion on $\cA^\times$. For $N > 0$, we apply 
\eqref{eq:leibniz} to $b = a^{-1}$ to obtain 
$$ 0 = \sum_S \dd\pi(x_S)(a) \dd\pi(x_{S^c})(a^{-1}), $$
where the sum is extended over all subsets of 
$\{1,\ldots, N\}$. This leads to 
$$ \dd\pi(x_1)\cdots \dd\pi(x_n)(a^{-1}) 
= -a^{-1} \sum_{ S \not=\eset} \dd\pi(x_S)(a) \dd\pi(x_{S^c})(a^{-1}). $$
That this sum depends, as an element of $\Mult^N(\g,\cA)$, continuously 
on $a$ follows by induction because $|S^c| < N$ whenever $S \not=\eset$. 
\end{prf}

\section{Smooth vectors for unitary representations} \label{sec:7}

In this section we provide a remarkably effective 
criterion for the smoothness of vectors in unitary representations 
of infinite dimensional Lie groups, namely that $v \in \cH^\infty$ 
if and only if the matrix coefficient $\pi^{v,v}(g) = \la \pi(g)v,v\ra$ 
is smooth in an identity neighborhood. 
This will be a simple consequence of the following 
observation.\begin{footnote}
{A finer analysis of the situation shows that 
if $K$ has continuous $k$-fold derivatives in both variables 
separately in some neighborhood of the diagonal, then 
$\gamma$ is a $C^k$-map, $k \in \N \cup \{\infty\}$ 
(cf.\ \cite[p.~78]{Kr49} for the case where $M$ is a real 
interval).}  
\end{footnote}

\begin{thm} \label{thm:kersmooth} Let $M$ be a smooth manifold, 
$\cH$ a Hilbert space and $\gamma \: M \to \cH$ be a map. 
Then $\gamma$ is a smooth map 
if and only if the kernel 
$K(x,y) := \la \gamma(x), \gamma(y) \ra$ is a smooth function 
on $M \times M$. 
\end{thm}

\begin{prf} Since the smoothness only refers to the real structure on 
$\cH$, we may assume that $\cH$ is a real Hilbert space. 
If $\gamma$ is smooth, then $K$ is also smooth because 
the scalar product on $\cH$ is real bilinear and continuous, hence smooth. 

{\bf Step 1:} $\gamma$ is continuous: This follows from 
$$ \|\gamma(x)- \gamma(y)\|^2 
= K(x,x) + K(y,y) - 2 K(x, y). $$

{\bf Step 2:} Now we consider the case where $M \subeq \R$ is an open 
interval, so that $\gamma \: M \to \cH$ is a curve in $\cH$. 
From Step 1 we know that 
$\gamma$ is continuous. For a fixed $t \in M$ we now consider on 
$M - t$ the function 
\begin{align*}
f(h) 
&:= \|\gamma(t+h)-\gamma(t)\|^2 
= K(t+h, t+h) + K(t,t) - 2 K(t+h,t). 
\end{align*} 
It is smooth with $f(0)= 0$ and 
$$ f'(0) = (\partial_1 K)(t,t) + (\partial_2 K)(t,t) 
- 2 (\partial_1 K)(t,t) = 0 $$
because $K$ is symmetric. 
This implies that 
$$ \lim_{h \to 0} \frac{2}{h^2} f(h) = f''(0) $$
exists, and from the Chain Rule and the symmetry of $K$ we obtain 
\begin{align*}
f''(0) 
&= (\partial_1^2 K)(t,t) + (\partial_2^2 K)(t,t) 
+ 2 (\partial_1 \partial_2 K)(t,t) 
- 2 (\partial_1^2 K)(t,t)\\ 
&= 2 (\partial_1 \partial_2 K)(t,t). 
\end{align*}
We conclude that 
\begin{equation}
  \label{eq:normlim}
\lim_{h \to 0}\Big\|\frac{1}{h}(\gamma(t+h)-\gamma(t))\Big\|^2 
= \lim_{h \to 0} \frac{1}{h^2} f(h) = \frac{1}{2} f''(0) 
=  (\partial_1 \partial_2 K)(t,t) 
\end{equation}
exists. We also note that, for each $s \in M$, 
$$ \frac{d}{dt} \la \gamma(t), \gamma(s)\ra 
= (\partial_1 K)(t,s) $$ 
exists. Since $\frac{1}{h}(\gamma(t+h)-\gamma(t))$ is bounded 
by \eqref{eq:normlim}, it follows that 
$ \frac{d}{dt} \la \gamma(t), v \ra$ 
exists for each $v$ in the closed subspace generated by $\gamma(M)$. 
For $v \in \gamma(M)^\bot$, the expression vanishes anyway, so that 
$\gamma'(t)$ exists weakly and satisfies 
$$ \la \gamma'(t), \gamma(s) \ra = (\partial_1 K)(t,s). $$
In particular, we have 
\begin{align*}
\|\gamma'(t)\|^2 
&= \lim_{h \to 0} \la \gamma'(t), \frac{\gamma(t+h)-\gamma(t)}{h}\ra 
= \lim_{h \to 0} \frac{1}{h}\big( (\partial_1 K)(t,t + h) - 
(\partial_1 K)(t,t)\big) \\
&= (\partial_1 \partial_2 K)(t,t) 
= \lim_{h \to 0}\Big\|\frac{1}{h}(\gamma(t+h)-\gamma(t))\Big\|^2, 
\end{align*}
and this implies that 
$$ \lim_{h \to 0} \frac{1}{h}(\gamma(t+h)-\gamma(t)) = \gamma'(t) $$
holds in the norm topology of $\cH$. 
We conclude that $\gamma$ is a $C^1$-curve.

Next we observe that 
$$ \la \gamma'(t), \gamma'(s) \ra 
= \partial_1 \partial_2 \la \gamma(t), \gamma(s) \ra 
= (\partial_1 \partial_2 K)(t,s)  $$
is also a smooth kernel. Therefore the argument above implies 
that $\gamma'$ is $C^1$, so that $\gamma$ is $C^2$. Repeating this 
argument shows that $\gamma$ is $C^k$ for every $k \in \N$, hence smooth.

{\bf Step 3:} Now we consider a general locally convex manifold $M$.  
As the assertion of the proposition is local, we may 
assume that $M$ is an open subset of a locally convex space. 
From Step 2 we derive that the map $\gamma$ has directional 
derivatives in all directions, which leads to a ``tangent map'' 
$$ \dd\gamma \: TM \to \cH, \quad 
(\dd\gamma)(x,v) := \derat0 \gamma(x + tv). $$ 
Since the kernel 
\begin{align*}
&\ \la (\dd\gamma)(x,v), (\dd\gamma)(y,w) \ra 
= \Big\la \derat0 \gamma(x + tv), \deras0 \gamma(y + sw) \Big\ra\\
&= \derat0 \deras0 \la \gamma(x + tv), \gamma(y + sw) \ra\\
&= \derat0 \deras0 K(x + tv, y + sw) 
= (\partial_{(v,0)} \partial_{(0,w)} K)(x,y),  
\end{align*}
is smooth, Step 1 shows that $\dd\gamma$ is continuous, so that 
$\gamma$ is a $C^1$-map. Applying the same argument to 
$\dd\gamma$ instead of $\gamma$, it follows that 
$\dd\gamma$ is $C^1$ and hence that $\gamma$ is $C^2$. Iterating 
this argument implies that $\gamma$ is smooth. 
\end{prf}

The following theorem substantially sharpens the well-known 
criterion (for finite dimensional groups) 
that $v \in \cH^\infty$ if its orbit map $\pi^v$ is weakly 
smooth, i.e., all matrix coefficients 
$\pi^{v,w}(g) = \la \pi(g)v,w\ra$ are smooth 
(cf.\ \cite[Cor.~10.1.3]{Sch90}). In \cite[p.~278]{Po74} 
one finds a remark suggesting its validity for finite dimensional 
groups, which is proved in \cite[Prop.~X.6.4]{Ne00} by 
using Goodman's characterization of smooth vectors 
(\cite{Go69}). 

\begin{thm} \label{thm:unitautsmooth} 
If $(\pi, \cH)$ is a unitary representation of a Lie group 
$G$, then $v \in \cH$ is a smooth vector if and only if the 
corresponding matrix coefficient $\pi^{v,v}(g) := \la \pi(g)v,v\ra$ 
is smooth on a $\1$-neighborhood in $G$. 
\end{thm}

\begin{prf} Clearly, $\pi^{v,v}$ is smooth if $v$ is a smooth vector. 
Suppose, conversely, that $\pi^{v,v}$ is smooth in a $\1$-neighborhood 
$U \subeq G$ and let $U'$ be a $\1$-neighborhood with 
$h^{-1}g \in U$ for $g,h \in U'$. 
In view of Theorem~\ref{thm:kersmooth}, the smoothness of $\pi^v$ 
on $xU'$, $x \in G$, is equivalent to the 
smoothness of the function 
$$ (g,h) \mapsto \la \pi(xg)v, \pi(xh)v \ra = \pi^{v,v}(h^{-1}g) $$ 
on $U' \times U'$, which follows from the smoothness of 
$\pi^{v,v}$ on $U$. 
\end{prf}

\begin{cor} \label{cor:a.14} If the continuous unitary representation 
$(\pi, \cH)$ of the Lie group $G$ has a cyclic vector $v$ for which the 
function $\pi^{v,v}(g) := \la \pi(g)v,v \ra$ is smooth on some identity 
neighborhood, then the representation $(\pi, \cH)$ is smooth, i.e., 
$\cH^\infty$ is dense. 
\end{cor}

\begin{prf} The preceding theorem implies that $v$ is a smooth 
vector, and therefore $\Spann \pi(G)v$ consists of smooth 
vectors. Hence $\cH^\infty$ is dense. 
\end{prf}

\begin{cor} \label{cor:7.4} 
If $\phi$ is a positive definite function on a Lie group 
$G$ which is smooth in a $\1$-neighborhood, then $\phi$ is smooth. 
\end{cor}

\begin{prf} Since $\phi$ is positive definite, the GNS construction 
provides a unitary representation $(\pi, \cH)$ of $G$ 
and a vector $v \in \cH$ with 
$\phi = \pi^{v,v}$.  Now Theorem~\ref{thm:unitautsmooth} implies that 
$v \in \cH^\infty$, but this implies that $\phi$ is smooth on all 
of $G$. 
\end{prf}

\section{$C^1$-vectors for Banach representations} \label{sec:9}

In this section we consider a continuous representation 
$(\pi, V)$ of the Banach--Lie group $G$ on the 
Banach space $V$, which implies in particular that 
$\pi$ is locally bounded (Lemma~\ref{lem:4.1}). 
Since Banach--Lie groups are locally exponential, 
Lemma~\ref{lem:a.11c} implies $v \in V$ is a $C^1$-vector if and only if 
$v \in \cD_\g$ and the map 
$\omega_v \: \g \to V, x \mapsto \oline{\dd\pi}(x)$ is continuous, 
which implies in particular that it is linear. 
The goal of this section is to see that 
the space $\cD_\g$ coincides with the space of $C^1$-vectors 
for the action of $G$ on $V$ (Theorem~\ref{thm:contlin}). 
In view of what we know already, the main point is that, for every 
$v\in \cD_\g$, the map $\omega_v$ is continuous. 
Surprisingly (for us), the most difficult part in our argument is 
to see that the maps $\omega_v \: \g \to V, x \mapsto 
\oline{\dd\pi}(x)v$, are additive for each $v \in \cD_\g$.

\begin{lem} \label{lem:1.1} Let $F \: [0,\infty[ \to \cL(V)$ be a curve 
with $F(0) = \1$ and $M > 0$ with 
$$  \|F(t)\| \leq M e^{at} \quad \mbox{ for } 
\quad t \geq 0 \quad \mbox{ and } \quad 
\|F(t)^n\| \leq M \quad \mbox{ for } \quad tn \leq 1. $$
Then there exists an $\omega \in \R$ with 
$$ \|F(t)^n\| \leq M e^{\omega nt} \quad\mbox{ for } \quad t\geq 0, n\in \N. $$
\end{lem}

\begin{prf} We put $\omega := 2\max(a, \log M)$ and observe that 
$M \geq \|F(0)\| = 1$ implies $\omega \geq 0$. 
The assertion clearly holds for $t = 0$, so that we may w.l.o.g.\ assume 
that $t > 0$. 

First we discuss the case where $t \leq 1$. 
Then we define $k \in \N$ by 
$\frac{1}{2^k} <  t \leq \frac{1}{2^{k-1}}$ and write 
$n$ as $n = q 2^{k-1} + r$ with $q,r \in \N_0$ and $r < 2^{k-1}$. 
Now $rt < 2^{k-1} t \leq 1$ leads to 
\begin{align*}
\|F(t)^n\| 
&\leq \|F(t)^{2^{k-1}}\|^q \|F(t)^r\| 
\leq M^q \cdot M. 
\end{align*}
On the other hand, 
$q \leq \frac{n}{2^{k-1}} < 2nt,$
so that 
$$ \|F(t)^n\| \leq M e^{q \log M} 
\leq M e^{nt 2 \log M}  \leq M e^{nt\omega}. $$

Now we consider the case where $t \geq 1$ and observe that 
\begin{align*}
 \|F(t)^n\| 
&\leq M^n e^{nat} 
= e^{n \log M + nat} \leq e^{n t \log M + nt a}
= e^{nt(\log M + a)} \leq e^{nt\omega} \leq M e^{nt\omega}. 
\end{align*}
\end{prf}

\begin{thm} \label{thm:additive} Let $(\pi, V)$ be a representation 
of the Banach--Lie group $G$ on the Banach space $V$ for which 
the corresponding action is continuous and $v \in \cD_\g$. 
Then the map 
$$ \omega_v \: \g \to V, \quad v \mapsto \oline{\dd\pi}(x)v $$
is linear. 
\end{thm}

\begin{prf} {\bf Step 1:} 
Since $G$ is a Banach--Lie group, 
the Trotter Product Formula 
$$ \lim_{n \to \infty} \big(\exp_G(tx/n)\exp_G(ty/n)\big)^n 
= \exp_G(t(x+y)) $$
holds uniformly for $|t| \leq N$ and any $N \in \N$. 

Below we need a slight refinement which can be obtained as follows. 
First we observe that in a Banach--Lie algebra, 
we have for each real sequence $\alpha_n \to s > 0$ and the 
Baker--Campbell--Hausdorff multiplication $*$ the relation 
$$ n\Big(\frac{\alpha_n}{n} x *\frac{\alpha_n}{n} y\Big) 
= \alpha_n 
 \frac{n}{\alpha_n} \Big(\frac{\alpha_n}{n} x *\frac{\alpha_n}{n} y\Big) 
\to s(x+y), $$
so that applying the exponential function on both sides leads to 
\begin{equation}
  \label{eq:trott-mod}
\Big(\exp_G(\alpha_n x/n) \exp_G(\alpha_n y/n)\Big)^n \to 
\exp_G(s(x+y)). 
\end{equation}

{\bf Step 2:} For $x,y \in \g$, we consider the strongly continuous curve 
$$ F \: \R \to \GL(V), \quad F(t) := \pi(\exp_G tx) \pi(\exp_G ty). $$ 
Since $\pi$ is locally bounded, it is bounded 
in a neighborhood of $\exp_G([0,1](x+y))$, so that we find with 
Step $1$ a $M_0 > 0$ with 
$$ \|F(t/n)^n\| \leq M_0 \quad \mbox{ for } \quad 0 \leq t \leq 1. $$
The strongly continuous one-parameter semigroups 
$\pi(\exp_G(tx))$ and $\pi(\exp_G(ty))$ satisfy 
$$ \|\pi(\exp_G(tx))\| \leq M_1 e^{t \omega_1}, \quad 
 \|\pi(\exp_G(ty))\| \leq M_2 e^{t \omega_2} $$
for suitable $M_1, M_2 > 0$, $\omega_1, \omega_2 \in \R$ 
and all $t \geq 0$ (\cite[Thm.~2.2]{Pa83}). This leads to the estimate 
$$ \|F(t)\| \leq M_1 M_2 e^{t(\omega_1 + \omega_2)}. $$
Applying Lemma~\ref{lem:1.1}, we now find a constant 
$M > 0$ and $\omega \in \R$ with 
\begin{equation}
  \label{eq:festi}
\|F(t)^n\| \leq M e^{\omega nt} \quad \mbox{ for } \quad 
t\geq 0, n \in \N. 
\end{equation}

{\bf Step 3:} From $\pi(g) \cD_x = \cD_{\Ad(g)x}$ for $g \in G$ and 
$x \in \g$ it follows that 
$\cD_\g$ is a $G$-invariant subspace of $V$. 
Since we are only interested in elements of $\cD_\g$, we may w.l.o.g.\ 
assume that $\cD_\g$ is dense in $V$. 
From the definition it immediately follows that 
$\oline{\dd\pi}(\lambda x)v = \lambda \oline{\dd\pi}(x)v$ 
for $\lambda \in \R$. 

For the operators $C := \oline{\dd\pi}(x)$ and 
$D := \oline{\dd\pi}(y)$ we then find that 
$\cD_\g \subeq \cD(C) \cap \cD(D)$ is dense in $V$. 
Further, \eqref{eq:festi} is the estimate required in 
\cite[Cor.~5.2(ii)]{Nek08}, which asserts that the unbounded operator 
$C + D$, defined on $\cD(C) \cap \cD(D)$, 
is closable and its closure generates 
a $C_0$-semigroup if and only if there exists a dense linear 
subspace $A \subeq \cD(C) \cap \cD(D)$ such that for all $f \in A$ 
and $s \geq 0$ there exist relatively compact 
sequences $(f_n^s)_{n \in \N}$ in $\cD(C) \cap \cD(D)$ such that 
\begin{itemize}
\item[\rm(a)] The set $\{ (C + D)f_n^s \: n \in \N\}$ is precompact 
for every $s \geq 0$. 
\item[\rm(b)] $\lim_{n \to \infty} \|F(t/n)^{[ns]}f  - f_n^s\| = 0$. 
\end{itemize} 
We put $A := \cD_\g$, and for $f \in \cD_\g$, we put 
$f_n^s := \pi(\exp_G(ts(x+y)))f$. Then (a) is trivially satisfied. 
To verify (b), we note that for $\alpha_n := \frac{[ns]}{n}$ the relation 
$ns - 1 \leq [ns] \leq ns$ implies 
$\alpha_n = [ns]/n \to s$. Hence 
\eqref{eq:trott-mod} leads to 
$$ F\big(\frac{\alpha_n t}{n}\big)^n f\to \pi(\exp_G(st(x+y)))f 
 \quad \mbox{ for } \quad  t \in \R, f \in \cD_\g. $$
We thus obtain 
$$ \lim_{n \to \infty} F(t/n)^{[ns]}f
= \lim_{n \to \infty} F(\alpha_n t/[ns])^{[ns]}f 
= \pi(\exp_G(st(x+y)))f = f_n^s. $$
Now \cite[Cor.~5.2]{Nek08} (see also the concluding Remark in 
\cite{Nek08b}) applies and yields 
$$ \oline{\oline{\dd\pi}(x) + \oline{\dd\pi}(y)} 
= \oline{C + D}  = \oline{\dd\pi}(x+y). $$
In particular, we obtain 
$\omega_v(x+y) = \omega_v(x) + \omega_v(y)$ for $v \in \cD_\g$.
\end{prf}

\begin{rem} \label{rem:a.10} 
If $\pi(G)$ consists of isometries, there are more direct 
arguments for the additivity of $\omega_v$. 
For the strongly continuous curve 
$$ F \: \R \to \GL(V), \quad F(t) := \pi(\exp_G tx) \pi(\exp_G ty), $$ 
we obtain for each $v \in \cD_\g$ the relation 
\begin{align*}
\frac{1}{t}(F(t)v-v) 
&= \pi(\exp_G tx)\frac{1}{t}\big(\pi(\exp_G(ty))v-v\big) 
+ \frac{1}{t}\big(\pi(\exp_G(tx))v-v\big) \\
&\to \oline{\dd\pi}(y)v + \oline{\dd\pi}(x)v,  
\end{align*}
so that the unbounded operator $F'(0)$ is defined on $\cD_\g$, 
where it coincides with 
$\oline{\dd\pi}(x)+ \oline{\dd\pi}(y)$. 

On the other hand, the Trotter--Product Formula in $G$ yields 
$$ \lim_{n \to \infty} F(t/n)^n 
=  \lim_{n \to \infty} \pi\Big( \big(\exp_G(x/t) \exp_G(y/t)\big)^n \Big)
=  \pi\big(\exp_G(t(x+y))\big) $$ 
in the strong operator topology. 
Now \cite[Theorem~3.1]{Chr74} implies that 
$\oline{\dd\pi}(x+y)$ extends the operator $F'(0)$ on 
$\cD_\g$. In particular, we obtain for $v \in \cD_\g$ that  
$\oline{\dd\pi}(x+y)v = \oline{\dd\pi}(x)v + \oline{\dd\pi}(y)v.$
\end{rem}

We now take a closer look at the continuity of the linear 
maps $\omega_v$ from Theorem~\ref{thm:additive}. 

\begin{lem} \label{lem:autocont} 
Let $(\pi, V)$ be a continuous representation 
of the Fr\'echet--Lie group $G$ on the Fr\'echet space $V$. 
Then, for each $v \in \cD_\g$ for which $\omega_v$ 
is linear, it is continuous. 
\end{lem}

\begin{prf} (S.~Merigon) Assume that $\omega_v$ is a linear map. 
In view of the Closed Graph Theorem (\cite[Thm.~2.15]{Ru73}), 
it suffices to show that the graph of $\omega_v$ is closed. 
Suppose that $x_n \to x$ in $\g$ such that 
$\omega_v(x_n) \to w$. Then we obtain from 
\eqref{eq:diff-rel0} the relation 
\begin{align*}
\pi(\exp_G(tx_n))v - v 
&= t\int_0^1 \pi(\exp_G(stx_n))\omega_v(x_n)\, ds. 
\end{align*} 
Since the function 
$$ [0,1]^2 \times \g \times V \to V, \quad 
(s,t,x,v) \mapsto \pi(\exp_G(stx))v $$
is continuous, integration over $s \in [0,1]$ leads to a continuous function 
$$ [0,1] \times \g \times V \to V, \quad 
(t,x,v) \mapsto \int_0^1 \pi(\exp_G(stx))v\, ds. $$
We thus obtain in the limit $n \to \infty$: 
$$ \pi(\exp_G(tx))v - v = t\int_0^1 \pi(\exp_G(stx))w\, ds. $$
For the derivative in $t = 0$ this leads to 
$$ \omega_v(x) = \oline{\dd\pi}(x)v = \int_0^1 \pi(\exp_G(0 x))w\, ds 
= \int_0^1 w\, ds = w. $$
\end{prf}

For continuous actions of Banach--Lie groups 
on Banach spaces we thus obtain with 
Theorem~\ref{thm:additive} and Lemmas~\ref{lem:autocont} and \ref{lem:a.11c}: 

\begin{thm} \label{thm:contlin} Let $(\pi, V)$ be a representation 
of the Banach--Lie group $G$ on the Banach space $V$ for which 
the corresponding action is continuous. 
Then $\cD_\g$ coincides with the space of $C^1$-vectors. 
\end{thm}

\section{$C^k$-vectors for Banach representations} \label{sec:10}

We continue our discussion of differentiable vectors for 
a continuous representation $(\pi,V)$ 
of a Banach--Lie group $G$ on a Banach space $V$. 
We know already that for each $v \in \cD_\g$ the map 
$\omega_v \: \g \to V$ is continuous and linear 
(Theorems~\ref{thm:additive}, Lemma~\ref{lem:autocont}).  
We thus obtain a norm on $\cD_\g$ by 
$$ \|v\|_1 :=  \|v\| + \|\omega_v\| = 
\|v\| + \sup_{\|x\| \leq 1} \|\omega_v(x)\| $$ 
(cf.\ \cite{Go70} and \cite{JM84} 
for similar constructions for finite dimensional 
Lie algebras). 
With respect to this norm on $\cD_\g$, the bilinear map 
$$ \g \times \cD_\g \to V, \quad 
(x,v) \mapsto \omega_v(x) = \oline{\dd\pi}(x)v $$
satisfies 
$\|\omega_v(x)\| \leq \|\omega_v\| \|x\| 
\leq \|v\|_1 \|x\|,$
so that it is continuous. 

For the sake of completeness, we recall the following 
variant of \cite[Lemmas~A.1/2]{Ne02}: 

\begin{lem} \label{lem:contfactor} Let 
$X_1, \ldots, X_n$, $Y$ and $Z$ be Banach spaces and 
$\eta \: Y \to Z$ a continuous injection. 
Suppose that $A \: X_1 \times \cdots \times X_n \to Z$ 
is a continuous $n$-linear map with 
$\im(A) \subeq \im(\eta)$. 
Then the induced $n$-linear map 
$$ \tilde A \: X_1 \times \cdots \times X_n \to Y 
\quad \hbox{ with } \quad \eta \circ \tilde A = A $$
is continuous. 
\end{lem}

\begin{prf} We argue by induction on $n$. 
First we consider the case $n = 1$. By the Closed Graph Theorem, 
it suffices to show that the graph of $\tilde A \: X_1 \to Y$ is closed. 
Assume that 
$$(x_n, \tilde A x_n) \to (x,y) \in X  \times Y.$$ 
Then $\eta(\tilde A x_n) = A x_n \to Ax$ 
implies that $\eta(y) = Ax = \eta(\tilde Ax)$, and therefore 
$\tilde A x = y$. Thus $\tilde A$ is continuous. 

Now we consider the general case $n > 1$. 
From the case $n = 1$,  we know that for 
fixed elements $x_j \in X_j$, $j < n$, the map 
$$X_n \to Y, \quad x \mapsto \tilde A(x_1,\cdots, x_{n-1},x) $$
is continuous. By the induction hypothesis, applied to the continuous 
inclusion 
$$\cL(X_n,Y) \to \cL(X_n,Z), \quad B \mapsto \eta \circ B,$$ 
the $(n-1)$-linear map 
$$ \tilde C \: X_1 \times \cdots \times X_{n-1} \to \cL(X_n, Y), \quad 
(x_1, \cdots, x_{n-1}) \mapsto 
\tilde A(x_1,\ldots, x_{n-1}, \cdot) $$
is continuous because the corresponding map 
$$ C \: X_1 \times \cdots \times X_{n-1} \to \cL(X_n, Z), \quad 
(x_1, \cdots, x_{n-1}) \mapsto A(x_1,\ldots, x_{n-1}, \cdot) $$
is continuous. Since the bilinear evaluation map 
$$ \cL(X_n, Z) \times X_n \to Z, \quad (\phi,x) \mapsto \phi(x) $$
is continuous, it follows that 
$\tilde A$ is a continuous $n$-linear map. 
\end{prf}

\begin{lem} \label{lem:domcomp} $\cD_\g$ is complete 
with respect to $\|\cdot\|_1$. 
\end{lem}

\begin{prf} The map $\iota \: \cD_\g\to V \times \cL(\g,V), v 
\mapsto (v, \omega_v)$ is isometric with respect to 
$\|\cdot\|_1$ on $\cD_\g$ and the norm 
$\|(v,\alpha)\| := \|v\| + \|\alpha\|$ on the product 
Banach space $V \times \cL(\g,V)$. Therefore the completeness of 
$\cD_\g$ is equivalent to the closedness of the graph 
$\Gamma(\omega) = \im(\iota)$ of the linear map 
$\omega \: \cD_\g \to \cL(\g,V), v \mapsto \omega_v.$
We now show that $\Gamma(\omega)$ is closed. 

Assume that $(v_n, \omega_{v_n}) \to (v, \alpha)$. 
Then \eqref{eq:diff-rel0} yields for each $x \in \g$ the relation 
$$ \pi(\exp_G(tx))v_n - v_n 
= t\int_0^1 \pi(\exp_G(stx)) \omega_{v_n}(x)\, ds. $$
Passing to the limit on both sides yields
$$ \pi(\exp_G(tx))v - v 
= t\int_0^1 \pi(\exp_G(stx)) \alpha(x)\, ds. $$ 
This implies that $v \in \cD_x$ with 
$\oline{\dd\pi}(x)v = \alpha(x).$ 
We conclude that $v \in \cD_\g$ with $\alpha = \omega_v$, 
and hence that $\Gamma(\omega)$ is closed. 
\end{prf}

\begin{lem}
  \label{lem:contmulti}
For each $v \in \cD_\g^n$, the corresponding $n$-linear map
$$\omega_v^n \: \g^n \to V, \quad 
(x_1,\ldots, x_n) \mapsto \oline{\dd\pi}(x_1)
\oline{\dd\pi}(x_2) \cdots \oline{\dd\pi}(x_n)v $$
is continuous. 
\end{lem}

\begin{prf} First we observe that Theorem~\ref{thm:additive}  
implies that $\omega_v^n$ is $n$-linear. We argue by induction on $n$ 
to show that it is continuous. 
For $n = 1$ this follows from Lemma~\ref{lem:autocont}. 
Now we assume $n>1$ and that $\omega_v^{n-1} \: \g^{n-1} \to V$ 
is a continuous $(n-1)$-linear map. From $v \in \cD_\g^n$ we 
derive that $\im(\omega_v^{n-1}) \subeq \cD_\g$. 
Since $\cD_\g$ is a Banach space with a continuous injection into  
$V$ (Lemma~\ref{lem:domcomp}), 
Lemma~\ref{lem:contfactor} implies that the 
$(n-1)$-linear map 
$\omega_v^{n-1} \: \g^{n-1} \to \cD_\g$
is continuous. Since the bilinear evaluation map 
$\g \times \cD_\g \to V, (x,w) \mapsto \omega_w(x)$ is 
continuous, we see that 
$$\omega_v^n(x_1, \ldots, x_n) 
= \omega_{\omega_v^{n-1}(x_2,\ldots, x_n)}(x_1) $$
is a continuous $n$-linear map. 
\end{prf}

\begin{thm}\label{thm:10} For a continuous 
representation of the Banach--Lie group $G$ on the Banach space $V$,  
the space $\cD_\g^k$ coincides with the space $V^k$ 
of $C^k$-vectors. 
This space is complete with respect to the norm 
  \begin{equation}
    \label{eq:k-norm}
 \|v\|_k := \|v\| + \sum_{j = 1}^k \|\omega_v^j\|,  
  \end{equation}
hence a Banach space, and the bilinear map 
\[ \xi \: \g \times \cD_\g^{k+1} \to \cD_\g^k, \quad 
(x,v) \mapsto \oline{\dd\pi}(x)v  \] 
is continuous with respect to the norms $\|\cdot\|_j$ on 
$\cD_\g^j$, $j = k,k+1$. 
\end{thm}

\begin{prf} 
Combining the preceding lemma with the general Lemma~\ref{lem:a.12}, 
it follows that $V^k = \cD_\g^k$. To verify the completeness of $\cD_\g^k$, 
we argue by induction. For $k = 1$, this is Lemma~\ref{lem:domcomp}. 
Assume that $k > 1$. We have to show that the graph of the linear map 
$$ \omega \: \cD_\g^k \to \prod_{j = 1}^k \Mult^j(\g,V), 
\quad v \mapsto (\omega^1_v,\ldots, \omega^k_v) $$
is closed in $V \times \prod_{j = 1}^k \Mult^j(\g,V)$. 
Suppose that 
$$(v_n, \omega_{v_n}^1, \ldots,\omega_{v_n}^k) \to 
(v,\alpha_1, \ldots, \alpha_k) 
\in V \times \prod_{j = 1}^k \Mult^j(\g,V).$$
Our induction hypothesis implies that $v \in \cD_\g^{k-1}$ 
with $\alpha_j = \omega_v^j$ for $j \leq k-1$. 

For $x, x_2,\ldots, x_k \in \g$ we further have 
$$ \oline{\dd\pi(x)} \omega_{v_n}^{k-1}(x_2, \ldots, x_k)  
= \omega_{v_n}^k(x, x_2, \ldots, x_k) \to  \alpha(x, x_2, \ldots, x_k)$$ 
and $\omega_{v_n}^{k-1}(x_2, \ldots, x_k)  
\to \omega_{v}^{k-1}(x_2, \ldots, x_k)$. 
Since the graph of $\oline{\dd\pi}(x)$ is closed (apply 
Lemma~\ref{lem:domcomp} with $\g = \R$),  
we obtain 
$$ \omega_{v}^{k-1}(x_2, \ldots, x_k) \in \cD_x 
\quad \mbox{ and } \quad 
\oline{\dd\pi}(x)\omega_{v}^{k-1}(x_2, \ldots, x_k) 
= \alpha(x, x_2, \ldots, x_k). $$
This implies that $v \in \cD_\g^k$ with $\alpha = \omega_v^k$. 

To show that the bilinear map $\xi$ is continuous, we first note that, 
for $v \in \cD_\g^{k+1}$, we have for $j \leq k$ the estimate 
\[ \|\omega_{\xi(x,v)}^j\| 
= \|\omega_{\oline{\dd\pi}(x)v}^j\| 
\leq \|\omega_v^{j+1}\| \|x\|, \] 
and this implies that 
\[ \|\xi(x,v)\|_k = \sum_{j = 0}^k \|\omega_{\xi(x,v)}^j\| 
\leq \sum_{j = 1}^{k+1} \|\omega_v^{j}\| \|x\| 
\leq \|v\|_{k+1} \|x\|.\] 
Therefore $\xi$ is continuous. 
\end{prf}

\begin{rem} Note that we have for each $n$ an isometric embedding 
$$ \cD_\g^n \into \cD_\g^{n-1} \times \Mult^n(\g,V), $$
where the norm on the product space is 
$\|(v,\alpha)\| := \|v\|_{n-1} + \|\alpha\|$. 
\end{rem}

\begin{rem} \label{rem:10.6} 
The preceding discussion leads in particular to a Fr\'echet 
structure on $\cD_\g^\infty$, 
considered as a subspace of $\prod_{n = 0}^\infty \Mult^n(\g,V)$. 
It coincides with the one from Definition~\ref{def:4.1}, so that 
we obtain a second proof of Proposition~\ref{prop:5.2}.  
\end{rem}

Now that we know that each space $\cD_\g^k$ is complete, 
it is natural to ask for the extent to which the $G$-action 
on this space is continuous. 

\begin{prop}
  \mlabel{prop:contonckvect} 
The representation $\pi^k$ of $G$ on the Banach space $V^k$ 
of $C^k$-vectors has the following properties: 
\begin{description}
\item[\rm(i)] $\pi^k$ is locally bounded. 
\item[\rm(ii)] For an element $v \in V^k$, the following 
are equivalent: 
\begin{description}
\item[\rm(a)] The maps $G \to \Mult^j(\g,V), g \mapsto \pi(g) \circ \omega_v^j$ 
are continuous for $j \leq k$. 
\item[\rm(b)] The maps $G \to \Mult^j(\g,V), g \mapsto 
\pi(g) \circ \omega_v^j \circ (\Ad (g)^{-1})^{\times j}$ 
are continuous for $j \leq k$. 
\item[\rm(c)] $v \in (V^k)^0$, i.e., the orbit 
map $G \to V^k, g \mapsto \pi(g)v$ is continuous. 
\item[\rm(d)] The orbit map $\pi^v \: G \to V$ is Fr\'echet-$C^k$, 
i.e., $v \in FV^k$. 
\end{description}
\item[\rm(iii)] $V^{k+1} \subeq FV^k$. 
\item[\rm(iv)] The subspace $FV^k$ of $V^k$ is closed 
and the $G$-action on $FV^k$ is continuous. It is the maximal 
$G$-invariant subspace of $V^k$ for which this is the case. 
In particular, the action of $G$ on $V^k$ is continuous  
if and only if $V^k = FV^k$. 
\end{description} 
\end{prop}

\begin{prf} (i) The group $G$ acts naturally by continuous linear operators 
on $\Mult^n(\g,V)$ via 
$$g.\omega 
:= \pi(g) \circ \omega \circ (\Ad(g^{-1}) \times \cdots \times \Ad(g^{-1}))
= \pi(g) \circ \omega \circ (\Ad(g^{-1})^{\times n}, $$
and we have 
$$ \|g.\omega\| \leq \|\Ad(g^{-1})\|^n \|\pi(g)\| \|\omega\|,  $$ 
i.e., the corresponding representation is locally bounded. 
Since the topological embedding 
\begin{equation}
  \label{eq:emb}
\omega \: V^k \to V \times \prod_{j = 1}^k \Mult^j(\g,V), \quad 
v \mapsto (v, \omega_v^1, \ldots, \omega_v^k) 
\end{equation}
is $G$-equivariant, the local boundedness of $\pi^k$ follows. 

(ii) The equivalence of (a) and (c) follows from the embedding 
\eqref{eq:emb}. 

Further, the equivalence of (a) and (b) follows from the 
fact that the action of $G$ on each space 
$\Mult^j(\g,V)$ by $g*\omega := \omega \circ (\Ad(g)^{-1})^{\times j}$ 
defines a morphism of Banach--Lie groups 
$G \to \GL(\Mult^j(\g,V))$ (cf.\ \cite[Exer.~IV.6]{Ne04}).  

For $v \in \cD_\g$, the orbit map $\pi^v \:  G\to V$ is $C^1$ 
with 
$$ \partial_{g.x}(\pi^v)(g) := T_g(\pi^v)(g.x) 
= \pi(g) \oline{\dd\pi}(x)v = \pi(g) \omega_v(x),$$ 
and, by iterating this argument, we obtain for 
$v \in V^k$ and $j \leq k$: 
$$ \partial_{g.x_1} \cdots \partial_{g.x_j}(\pi^v)(g) 
= \pi(g) \omega_v^j(x_1, \ldots, x_j). $$

With similar arguments as in the proof of 
\cite[Thm.~I.7]{Ne01a}, we now see that the $C^k$-map 
$\pi^v$ is $C^k$ in the Fr\'echet sense if and only if, 
for $j \leq k$, the maps 
$$ G \to \Mult^j(\g,V), \quad 
g \mapsto \big( (x_1, \cdots, x_j)\mapsto 
(\partial_{g.x_1} \cdots \partial_{g.x_j}\pi^v)(g)\big) $$
are continuous. In view of the preceding calculations, 
this means that (a) is equivalent to (d). 

(iii) Follows from the general fact that each 
$C^{k+1}$-map is $C^k$ in the Fr\'echet sense 
(cf.~\cite[Thm.~I.7(ii)]{Ne01a}, \cite{GN10}). 

(iv) Since the $G$-representation on $V^k \subeq \prod_{j = 0}^k 
\Mult^j(\g,V)$ is locally bounded, 
the subspace  of elements with continuous orbit maps  
is closed. As we have seen in (iii) above, this subspace coincides 
with $FV^k$. It is clearly invariant and Lemma~\ref{lem:4.1} 
implies that $G$ acts continuously on $FV^k$. 
\end{prf}

\begin{rem} \label{rem:9.8} 
If $G$ is finite dimensional, then each $C^k$-map is also 
Fr\'echet-$C^k$, so that the action of $G$ 
on $\cD_\g^k$ is again a continuous action 
(Remark~\ref{rem:5.3}). 
In sharp contrast to this situation is the fact that, 
for an infinite dimensional Banach--Lie group, 
$\cD_\g$ need not contain any non-zero continuous 
vector (cf.\ Remark~\ref{rem:11.6} below).  
\end{rem}

\begin{ex}
For the examples discussed in Section~\ref{sec:11} below, 
it is shown that the action of $G$ on $\cD_\g$ 
is continuous for $p \geq 4$, because in this 
case $\cD_\g^2$ contains the dense subspace $L^\infty([0,1])$ 
(cf.\ Remarks~\ref{rem:5.3} and \ref{rem:11.6}). 
\end{ex}

\begin{ex} We consider the action of the one-dimensional 
Lie group $G = \R$ on the Banach space $V = C(\R,\R)_{\rm per}$ 
of $1$-periodic functions by $(\pi(g)f)(x) := f(g+x)$. 
Then the uniform continuity of each element of $V$ implies 
that the orbit maps are continuous, and since $G$ acts by 
isometries, Lemma~\ref{lem:4.1} implies the continuity of 
the $G$-action on $V$. 

Since the point evaluations are continuous, the space
$\cD_\g^k$ of $C^k$-vectors consists of $C^k$-functions on 
$\R$. If, conversely, $f$ is a $C^1$-function, then 
$$ \frac{f(x + h) - f(x)}{h} = \int_0^1 f'(x + sh)\, ds $$
converges uniformly to $f'(x)$ for $h \to 0$, so that 
$f$ is a $C^1$-vector. A similar argument, based on the 
integral form of the remainder term in the Taylor formula 
implies that $\cD_\g^k = C^k(\R,\R)_{\rm per}$. 
\end{ex}

\section{A family of interesting examples} \label{sec:11}

We take a closer look at the unitary representation 
of the Banach--Lie group $G := (L^p([0,1],\R),+)$, $p \in [1,\infty[$,  
on the Hilbert space $\cH = L^2([0,1],\C)$ by $\pi(g)f := e^{ig}f$. 
In \cite{BN08} it is shown that, for $p = 2$, this 
representation is continuous with $\cH^\infty = \{0\}$. 
Here we shall see that it is always continuous and determine 
the space of $C^k$-vectors. In the following 
we write $\g := L^p([0,1],\R)$ for the Lie algebra of $G$ 
and abbreviate $L^p([0,1]) := L^p([0,1],\C)$. 

We start with a general observation on the inclusions between 
$L^p$-spaces. 

\begin{rem}\label{rem:11.1}
(a) If $(X,\mu)$ is a finite measure space, then 
$$ L^q(X,\mu) \subeq L^p(X,\mu)  \quad \mbox{ for } \quad 1 \leq p < q, $$
where the inclusion is continuous. 
In fact, for any measurable function $f \: X \to \C$, we have 
$$ \int |f|^p 
= \int_{\{|f|\leq 1\}}|f|^p + \int_{\{|f| > 1\}}|f|^p
\leq \mu(X)  + \int_{\{|f| > 1\}}|f|^q. $$  
This implies that $\|\cdot\|_p$ is bounded on the unit ball of 
$L^q(X,\mu)$. 

(b) Assume that $X$ has a decomposition into a sequence 
$(X_n)$ of pairwise disjoint subsets with 
$0 < \mu(X_n) \leq 2^{-n} \mu(X)$. We claim that 
$$ L^q(X,\mu) \not= L^p(X,\mu)  \quad \mbox{ for } \quad 1 \leq p < q.$$

We consider the function $f$ whose value on $X_n$ is constant 
$\mu(X_n)^{-1/q}$. Then \break 
$\|f\|_q^q = \sum_n \mu(X_n)^{-1} \mu(X_n) = \infty$ and 
$$\|f\|_p^p 
= \sum_n \mu(X_n)^{1-p/q} \leq 
\mu(X)^{1-p/q} \sum_n (2^{p/q-1})^n < \infty.$$  

(c) If $X = [0,1]$ is the unit interval and $\mu$ is Lebesgue measure, 
then the property under (b) is satisfied for every 
subset $Y \subeq X$ of positive measure. 
\end{rem}

\begin{lem}\label{lem:11.2}  The representation 
$(\pi, \cH)$ is continuous for every 
$p \geq 1$. 
\end{lem}

\begin{prf}  (cf.\ \cite[Prop.~2.1]{BN08} for $p = 2$)  
Let $(g_n)$ be a sequence in $G = L^p([0,1],\R)$ converging 
to $0$.  Then there exists a subsequence 
$(g_{n_k})_{k \in \N}$ converging almost everywhere 
to $0$. Then $|e^{ig_{n_k}} -1|^2\to 0$ almost everywhere, so 
that the Dominated Convergence Theorem implies that 
$\int_{[0,1]} |e^{ig_{n_k}} -1|^2 \to 0.$
For every $\xi \in L^\infty([0,1],\C)$, we thus obtain 
$$ \|\pi(g_{n_k})\xi - \xi\|_2^2 
= \int_{[0,1]} |e^{ig_{n_k}}-1|^2 |\xi|^2 
\leq \|\xi\|_\infty^2  \int_{[0,1]} |e^{ig_{n_k}}-1|^2 \to 0.$$
We conclude that 
$\pi(g_{n_k}) \to \1$ with respect to the strong topology. 

If $\pi$ is not continuous, then there exists a sequence 
$h_n \to 0$ in $G$ with $\pi(h_n) \not\to \pi(0)$. 
This means that there exist a $\1$-neighborhood $U$ in $\U(\cH)$ 
with respect to the strong operator topology  
for which $\{n \: \pi(h_n) \not\in U\}$ is infinite. 
This leads to a sequence $g_n$ in $G$ with 
$g_n \to 0$ and $\pi(g_n) \not\in U$ for every $n$, 
and we thus obtain a contradiction to the argument in the 
preceding paragraph. We conclude that $\pi$ is continuous.
\end{prf}

Since the infinitesimal generator of the one-parameter group 
$t \mapsto \pi(tf)$ is the multiplication with the function 
$f$, it follows that 
$$ \cD_\g = \cD_\g^1 = \{ \xi \in \cH \: (\forall f \in L^p([0,1],\R))\, f\xi 
\in \cH\}. $$

\begin{lem} We have 
$$ \cD_\g = 
\begin{cases}
\{0\} \text{ for } p < 2 \\
L^\infty([0,1])  \text{ for } p = 2 \\
L^{\frac{2p}{p-2}}([0,1]) \text{ for } p > 2. 
\end{cases}$$
\end{lem}

\begin{prf} Let $\xi \in \cD_\g$. 
If $\xi\not=0$, then there exists an 
$\eps > 0$ for which $X_\eps := \{|\xi| > \eps\}$ has positive 
measure. Now 
$L^p([0,1]) \cdot \xi \subeq L^2([0,1])$ 
implies that 
$L^p(X_\eps) \subeq L^2(X_\eps)$, 
and in view of Remark~\ref{rem:11.1}(b/c), this leads to $p \geq 2$. 

For $p = 2$, the multiplication with $\xi$ defines a bounded 
operator on $L^2([0,1])$, so that $\xi \in L^\infty([0,1])$. 
This proves that $\cD_\g = L^\infty([0,1])$ in this case. 

Now we assume that $p > 2$. The condition 
$\xi \cdot L^p([0,1]) \subeq L^2([0,1])$ is equivalent to 
$\xi^2 \cdot L^{p/2}([0,1]) \subeq L^1([0,1])$, which 
is equivalent to $\xi^2 \in L^q([0,1])$ for 
$\frac{1}{q} + \frac{2}{p} = 1,$
i.e., $q = \frac{p}{p-2}.$ This proves that 
$\cD_\g = L^{\frac{2p}{p-2}}([0,1])$. 
\end{prf}

Next we note that 
$\Spann \{ f^k \: f \in \g \} 
= \Spann \{ f_1 \cdots f_k \: f_1, \ldots, f_k \in \g \}$ 
shows that 
$$ \cD_\g^k = \{ \xi \in \cH \: (\forall f \in L^p([0,1],\R))\, f^k\xi 
\in \cH\}. $$
Now 
$$ \{ f^k \: f \in \g = L^p([0,1])\} 
= L^{p/k}([0,1]) \quad \mbox{ for } \quad p \geq k, $$
leads for $p \geq k$ to 
\begin{equation}
  \label{eq:e.1}
\cD_\g^k = \{ \xi \in \cH \: (\forall f \in L^{p/k}([0,1],\R))\, f\xi 
\in \cH\} =  \cD_{L^{p/k}}. 
\end{equation}

\begin{prop} \label{prop:11.4}
For $p \in \N$ and $k \leq p$ we have 
$$\cD^{k}_\g = \cD_{L^{p/k}}
= \begin{cases}
\{0\} \text{ for } k > \frac{p}{2} \\
L^\infty([0,1])  \text{ for } k = \frac{p}{2} \\
L^{\frac{2p}{p-2k}}([0,1]) \text{ for } k < \frac{p}{2}. 
\end{cases}$$
\end{prop}

For $p = 2$ we obtain in particular 
$\cD^2_\g = \{0\}$ and $\cD_\g^1 = L^\infty([0,1])$ 
which refines the observations in \cite{BN08}. 

\begin{rem}
The preceding proposition also shows that there exists 
for every $n \in \N$ a continuous unitary representation 
$(\pi, \cH)$ of a Lie group $G = (L^{2n}([0,1],\R),+)$ 
with $\cD_\g^{n+1} = \{0\}$ and $\cD_\g^n \not=\{0\}$. 
\end{rem}

\begin{rem} It is easy to see that the norm 
$\|\cdot\|_k$ on the space $\cD_\g^k$ is equivalent to the natural 
norm suggested by Proposition~\ref{prop:11.4}. 

For $p = 4$ and $\cD_\g = L^4([0,1])$, the density of 
$\cD^2 = L^\infty([0,1])$ in $\cD_\g$ implies the continuity 
of the isometric action of $G$ on $\cD_\g$ 
(cf.\ Remark~\ref{rem:5.3}). 
\end{rem}

\begin{rem} \label{rem:11.6} 
For maps between Fr\'echet spaces, we also have 
the stronger notion of $C^k$-maps in the Fr\'echet 
sense. If $\xi \in \cH = L^2([0,1],\C)$ is a $C^1$-vector 
for $G = (L^p([0,1],\R),+)$, then the orbit map 
$\pi^\xi(g) = e^{ig}\xi$ has in $g$ the differential 
$$ T_g(\pi^\xi)f = e^{ig}f\xi. $$
Therefore $\xi$ is a Fr\'echet-$C^1$-vector if and only if the map 
$$ F \: G \to \cL(\g,\cH), \quad g \mapsto M_{e^{ig}\xi}, \quad 
M_h f = hf, $$
is continuous. 

We consider the case $p = 2$. 
Then $\cL(\g,\cH) \cong \cL(\cH)$, so that we may consider 
$F$ as a map 
$$ F \: G \to L^\infty([0,1],\C), \quad g \mapsto e^{ig}\xi, $$
and the question is when this map is continuous. 

First we consider the case $\xi=1$. Then 
$F$ is a homomorphism of Banach--Lie groups. If it is 
continuous, it is smooth, which implies that 
$\L(F) \: \g \to L^\infty([0,1],\C), g \mapsto ig,$ 
is a continuous liner map, which is not the case. 

For a general $\xi \in L^\infty([0,1],\C)$, we 
consider for $\eps > 0$ the subsets 
$X_\eps := \{|\xi| \geq \eps\}$. On each of these sets,  
$\xi\res_{X_\eps}$ is invertible in the Banach algebra 
$L^\infty([0,1],\C)$, so that multiplication with 
$\xi^{-1}$ leads to the situation of the previous paragraph. 
Therefore $X_\eps$ has measure zero, and since $\eps > 0$ 
was arbitrary, $\xi =0$. Hence all Fr\'echet-$C^1$-vectors 
for $G = L^2([0,1],\R)$ are trivial. Another way to put this is 
to say that the action of $G$ on the Banach space 
$\cD_\g = L^\infty([0,1],\C)$ has no non-zero continuous vector. 

On the other hand, we know that every $C^2$-vector is a 
Fr\'echet-$C^1$-vector (\cite[Thm.~I.7]{Ne01a}), so that we obtain non-trivial 
Fr\'echet-$C^1$-vectors for $G = L^4([0,1],\R)$.
\end{rem}

\section{Smooth vectors for direct limits} \label{sec:12}

\begin{defn} If $(G_n)_{n \in \N}$ is a sequence of finite dimensional 
Lie groups with homomorphisms $\phi_n \: G_n \to G_{n+1}$, 
then Gl\"ockner has shown in \cite{Gl03} that the corresponding direct limit 
group $G := \indlim G_n$, endowed with the direct limit 
topology carries a compatible Lie group structure with 
Lie algebra $\g := \indlim \L(G_n)$, endowed with the direct limit 
topology. We call $G$ a {\it direct limit Lie group}. 
\end{defn}

\begin{rem} The Lie groups obtained by this construction are precisely the 
Lie groups $G$ with a smooth exponential function 
whose Lie algebra $\g$ is a countable union of finite dimensional 
subalgebras, endowed with the direct limit topology. 
In fact, writing $\g = \bigcup_{n \in \N} \g_n$ with 
finite dimensional Lie algebras $\g_n \subeq \g_{n+1}$, we obtain a 
corresponding sequence of Lie groups $G_n$ injecting into $G$, 
and then it is not hard to verify that the corresponding morphism of 
Lie groups $\indlim G_n \to G$ is an isomorphism 
(cf.\ \cite{GN10}). 
\end{rem}

\begin{thm} For each continuous unitary representation 
$(\pi, \cH)$ of a direct limit Lie group $G$, 
the space of smooth vectors is dense. 
\end{thm} 

\begin{prf} Since $\cH$ is a direct sum of subspaces on which the representation if cyclic, we may w.l.o.g.\ assume that the representation is cyclic. 
Since smoothness of a vector in $\cH$ is equivalent to smoothness for the 
identity component, we may also assume that $G$ is connected. 
Then $G$ is a countable direct limit of connected finite dimensional 
Lie groups, hence separable, and therefore the cyclicity of 
$(\pi,\cH)$ implies that $\cH$ is separable.

Danilenko shows in \cite{Da96} that there exists a dense 
subspace $\cD \subeq \cH$ satisfying 
\begin{description}
\item[\rm(a)] $\cD$ is $\pi(G)$-invariant.  
\item[\rm(b)] $\cD \subeq \cD_\g$.
\item[\rm(c)] $\oline{\dd\pi}(x)\cD \subeq \cD$ for every $x \in \g$. 
\end{description}

Since each group $G_n$ is in particular a Banach--Lie group, 
conditions (b) and (c) imply that $\cD \subeq \cD^\infty_{\g_n}$, 
so that Lemma~\ref{lem:a.12} implies that $\cD$ consists of 
smooth vectors for each $G_n$. Since $G$ is also the direct limit 
of the $G_n$ in the category of smooth manifolds 
(\cite{Gl03}), $\cD \subeq \cH^\infty$ 
consists of smooth vectors for $G$. 
\end{prf}

\begin{rem} Typical examples of Lie algebras $\g$ which  
are locally finite in the sense that every finite subset generates 
a finite dimensional subalgebra are nilpotent Lie algebras. 
These Lie algebras are countable direct limits of finite dimensional  
ones if they are countably dimensional. This condition is very 
restrictive, so that one is also interested in situations, where 
the Lie algebra $\g$ is not of countable dimension but 
carries a locally convex topology which is coarser than the direct 
limit topology. 

An important class of corresponding groups are the Heisenberg 
groups $\Heis(V)$ of a locally convex space $V$, endowed with a 
continuous scalar product (a locally convex euclidean space). 
More precisely, we have 
$$ \Heis(V) = \R \times V \times V, \quad 
(z,v,w)(z',v',w') := 
(z + z' + \shalf(\la v,w'\ra - \la v',w\ra), v + v', w + w'). $$

A unitary representation $(\pi, \cH)$ of this group with 
$\pi(z,0,0) = e^{iz}$ for $z \in \R$ provides a representation of 
the {\it canonical commutation relations} in the sense that 
the restriction to the subgroups 
$\{0\} \times V \times \{0\}$ and $\{(0,0)\} \times V$ 
defines a unitary representation, and for $v,w \in V$ we have 
$$ \pi(0,v,0) \pi(0,0,w) = e^{i\la v,w\ra} \pi(0,v,w) 
= e^{2i\la v,w\ra}  \pi(0,0,w) \pi(0,v,0). $$

In \cite{Heg72} Hegerfeldt shows that, if $V$ is separable, barreled 
and nuclear, for any continuous 
unitary representation of $G = \Heis(V)$ there exists a dense subspace 
$\cD$ with the following properties: 
\begin{description}
\item[\rm(a)] $\cD$ is $\pi(G)$-invariant.  
\item[\rm(b)] $\cD \subeq \cD_\g$.
\item[\rm(c)] $\oline{\dd\pi}(x)\cD \subeq \cD$ for every $x \in \g$, 
in particular $\cD \subeq \cD_\g^\infty$. 
\item[\rm(d)] For each $v \in \cD$ the map 
$\omega_v^n \: \g^n \to \cH$ is continuous. 
\end{description}
In view of Lemma~\ref{lem:a.12} and the local exponentiality of $G$, 
(c) and (d) imply that 
$\cD$ consists of smooth vectors for $G$. 

Actually Hegerfeldt shows that the elements of $\cD_\g$ even have analyticity 
properties that lead to holomorphic extensions of their orbit 
maps to the complexified group $G_\C$. 
For a detailed discussion of this aspect we refer to 
\cite{Ne10a}. 
\end{rem}

\section{Smooth vectors for projective limits} \label{sec:13}

Structurally direct limits of finite dimensional Lie groups are 
groups with a relatively simple structure, but they have many 
continuous unitary representation because they carry a very fine 
topology. This situation is the opposite of what we find for 
projective limits of finite dimensional Lie groups. 
These groups are also called {\it pro-Lie groups}, 
and the Lie groups among the pro--Lie groups have been characterized 
recently in \cite{HN09}. For any such Lie group 
$G = \prolim G_j$, it makes sense to ask for smooth vectors 
in unitary representations. As we shall see in this section, 
for this class of Lie groups the space $\cH^\infty$ is always dense. 
Here the main point is a quite general argument concerning 
unitary representations of projective limits 
of topological groups. 

Let $G$ be a topological group and ${\cal N} := (N_i)_{i \in I}$ be a filter basis of 
closed normal subgroups $N_i \trile G$ with $\lim {\cal N} = \{\1\}$,
i.e., for each $\1$-neighborhood $U$ in $G$ there exists some $i\in I$ with 
$N_i \subeq U$. 

\begin{lem}
  \label{lem:7.1} 
Let $(\pi, {\cal H})$ be a continuous unitary representation 
of $G$. Then the union $\bigcup_{i \in I} {\cal H}^{N_i}$ of the closed invariant 
subspaces ${\cal H}^{N_i}$ of $N_i$-fixed vectors is dense in ${\cal H}$. 
\end{lem}

\begin{prf}
Let $v \in {\cal H}$ and pick $\eps >0$. 
Then there exists some $i$ with $\pi(N_i)v \subeq B_\eps(v)$. 
Then $\oline{\conv(N_i.v)}$ is a bounded closed convex $N_i$-invariant subset of 
${\cal H}$, hence contains 
a $N_i$-fixed point by the Bruhat--Tits Fixed Point Theorem 
(\cite{La99}). 
Therefore $\dist(v, {\cal H}^{N_i}) \leq \eps$. 
\end{prf}

\begin{thm}
  \label{thm:7.2}
Any continuous representation $(\pi, {\cal H})$ of $G$ 
is a direct sum of representations on which some $N_i$ acts trivially. 
\end{thm}

\begin{prf}
Let ${\cal F}$ be the set of all sets $F$ of pairwise orthogonal 
$G$-invariant closed subspaces of ${\cal H}$ on which some $N_i$ acts trivially. 
We order ${\cal F}$ by set inclusion. Then ${\cal F}$ is inductively ordered, 
so that Zorn's Lemma implies the existence of a maximal element $F_m$. 
Then 
$${\cal K} := \oline{\bigoplus_{{\cal H}_j \in F_m} {\cal H}_j}$$ 
is a closed $G$-invariant subspace of ${\cal H}$. 
We claim that ${\cal K} = {\cal H}$, which implies the assertion. 
If this is not the case, 
${\cal K}^\bot$ is non-zero, and Lemma~\ref{lem:7.1} implies that 
for some $i$ the set $({\cal K}^\bot)^{N_i}$ is non-zero, contradicting the 
maximality of $F_m$. 
\end{prf}

\begin{cor} Each irreducible continuous unitary representation 
of $G$ factors through some $G/N_i$, i.e., the set $\hat G$ of 
equivalence classes of irreducible continuous unitary representations 
satisfies $\hat G = \bigcup_{i \in I} (G/N_i)\,\hat{}.$
\end{cor} 

\begin{thm} If $G = \prolim G_j$ is a Lie group which, as a topological 
group, is a projective limit of finite dimensional Lie groups 
$G_j$, then for each continuous unitary representation 
$(\pi,\cH)$ of $G$ the space $\cH^\infty$ of smooth vectors is 
dense. 
\end{thm} 

\begin{prf} Let $q_j \: G \to G_j$ be the natural projections 
and apply Theorem~\ref{thm:7.2} to the family $N_j = \ker q_j$. 
This reduces the problem to the case where some $N_j$ acts 
trivially on $\cH$, so that we actually have a representation of the 
finite dimensional quotient Lie group $G_j \cong G/N_j$. Now the 
assertion follows from the density of smooth vectors 
for  $G_j$ in~$\cH^{N_j}$ (\cite{Ga47}). 
\end{prf}

\begin{rem} The group $G= \R^\N$ is a projective limit 
of the Lie group $G_n = \R^n$, where $q_n \: G \to G_n$ 
is the projection onto the first $n$ factors. 
In this sense Example~\ref{ex:4.7} is a continuous unitary 
representation of a pro-Lie group. 
\end{rem}

\begin{rem}
 Let $G = (V,+)$ be the additive group of a locally convex space $V$. 
For each continuous seminorm $p \in \cP(V)$, we have a 
closed subspace $N_p := p^{-1}(0)$ for which $p$ induces a norm 
on the quotient space $V/N_p$. Now 
$\cN = \{ N_p \:p \in \cP(V)\}$ is a filter basis of closed subgroups 
with $\lim \cN = \{\1\}$, so that Theorem~\ref{thm:7.2} applies. 
We conclude that every continuous unitary representation 
$(\pi,\cH)$ of $V$ is a direct sum of representations 
$(\pi,\cH_i)$ on which some $N_p$ acts trivially, so that the 
representation $\pi_i$ factors through a representation 
of the normed space $V/N_p$. 
\end{rem}

\section{Perspectives} \label{sec:14}

There are several interesting problems concerning representations 
of infinite dimensional Lie groups $G$ on Banach spaces $V$. 

\begin{prob} (Integrability) Suppose that 
$V$ is a locally convex space, $\cD \subeq V$ a dense subspace 
and $\rho \: \g \to \End(\cD)$ a representation of the 
Lie algebra $\g$ on $\cD$. 

We thus obtain for each $x \in \g$ an unbounded operator 
$\rho(x)$ on $V$. We assume that all these operators are closable 
and that each closure $\oline{\rho(x)}$ generates a strongly 
continuous one-parameter group. A characterization of such 
operators on Banach spaces is given by the Hille--Yoshida Theorem, 
and \cite{Ko68} contains some generalization to locally convex spaces. 

Suppose that $\g$ is a Banach--Lie algebra. Then we obtain a map 
$$ F \: \g \to \GL(V), \quad x \mapsto e^{\oline{\rho(x)}}. $$
When does this map define a local group homomorphism, hence a representation 
of any corresponding simply connected Lie group $G$? More precisely, 
under which conditions on $\rho$ do we have 
$$ F(x*y) = F(x) F(y) $$
for $x$ and $y$ in a small ball centered in~$0$? 

For finite dimensional Lie algebras, problems of this kind 
are studied in \cite[Thms.~1.1]{Jo83} and Chapter~$8$ of \cite{JM84}. 
Maybe some of these results, such as Theorem~A.1-3 in \cite{Jo88} 
can be extended to Banach--Lie algebras. 
It would also be interesting to have a version of 
the integrability result \cite[Thm.~8.1]{JM84} for representations 
on locally convex space or \cite[Thm.~8.6]{JM84} for representations 
on Banach spaces (cf.\ also \cite{Mo65}). 
A first step in this direction is taken by S.~Merigon in 
\cite{Me10}, where he obtains such a result for representation 
on Hilbert spaces. 
\end{prob} 

\begin{prob} (Smoothness) We have seen above that for every continuous 
representation $(\pi,V)$ of the Banach--Lie group $G$ on 
the Banach space $V$, we obtain a sequence 
$(V^k, \|\cdot\|_k)$ of Banach spaces, where 
$V^k = \cD_\g^k$ is the space of $C^k$-vectors, endowed with its natural 
norm (Theorem~\ref{thm:10}), and in this picture $V^\infty$ is the 
projective limit of the Banach spaces $V^k$ 
(in the category of locally convex spaces) (cf.\ \cite[p.~20]{JM84} 
for finite dimensional Lie algebras). 

The $C^k$-variant of the derived action is given by the sequence 
\begin{equation}
  \label{eq:bilseq}
\g \times V^k \to V^{k-1},\quad (x,v) \mapsto \oline{\dd\pi}(x)v, \quad k \in \N, 
\end{equation}
of continuous bilinear maps. Interesting questions in this context are: 
\begin{description}
\item[\rm(a)] Does the density of $V^\infty$ in $V$ imply the density 
in each of the Banach spaces $V^k$? 
\item[\rm(b)] In \cite{JM84} the continuity of the action of 
$G$ on the Banach space $(\cD_\g, \|\cdot\|_1)$ (graph density) 
plays an important role. As follows from 
Proposition~\ref{prop:contonckvect} and Remark~\ref{rem:5.3}, 
this is equivalent to the density of the continuous vectors in 
$\cD_\g$, which in turn follows from the density of $\cD_\g^2$ in 
$\cD_\g^1$. Maybe these conditions can also be exploited for 
Banach--Lie groups. 
\item[\rm(c)] Suppose that we are given a Lie algebra representation $\rho \: \g \to \End(\cD)$, where 
$\cD \subeq V$ is a dense subspace. Suppose that all the 
$k$-linear maps $\omega_v^k \: \g^k \to V$ are continuous and write 
$V^k$ for the completion of $\cD$ with respect to the norm 
$$ \|v\|_k := \|v\| + \sum_{j = 1}^k \|\omega_v^j\|. $$
Then $V^k$ can be realized as subspaces of $V$. 
Is it possible to characterize in this context representations which 
are integrable to continuous representations of $G$ on $V$ 
(cf.\ \cite{JM84} for a discussion of similar problems for finite 
dimensional Lie algebras). Natural assumptions in this context are 
that that the closures of the operators $\rho(x)$ generate one-parameter 
groups preserving $\cD$ (and all the spaces $V^k$). 
\end{description}
\end{prob}

\begin{prob} We have seen in Sections~\ref{sec:4} and \ref{sec:5} 
that for every 
unitary representation $(\pi,\cH)$ of a Banach--Lie group 
$G$, the space $\cH^\infty$ carries a natural Fr\'echet topology with 
respect to which $G$ acts smoothly. Can this information be used to 
show in certain situations that a direct integral decomposition 
$$ (\pi, \cH) = \int_X^{\oplus} (\pi_x, \cH_x)\, d\mu(x) $$
also yields a ``direct integral decomposition'' 
$$ \cH^\infty = \int_X^{\oplus} \cH_x^\infty\, d\mu(x)? $$
This would be extremely useful for the analysis of smooth 
unitary representations. 
For results of this type for finite dimensional groups we refer 
to \cite{Ar76}. 
\end{prob}

\begin{prob} The argument in the proof of Theorem~\ref{thm:additive} 
touches on an interesting question concerning the differentiability 
of functions $f \: G \to \R$ on a Banach--Lie group. 
Suppose that for every $g \in G$ and $x \in \g$ the derivative
$$ \dd f(g)(g.x) := \derat0 f(g \exp_G(tx)) $$
exists. When are the maps $\dd f(g) \: \g \to \R$ linear? 
They clearly satisfy $\dd f(g)(\lambda x) = \lambda \dd f(g)x$ for 
$\lambda \in \R$, so that the additivity is the crucial 
issue. 

The connection to Theorem~\ref{thm:additive} is given by 
functions of the form 
$f(g) = \alpha(\pi(g)v)$ with $\alpha \in V'$ and $v \in \cD_\g$, 
because in this case we have 
$\dd f(\1)(x) = \alpha(\oline{\dd\pi}(x)v) = \alpha(\omega_v(x))$, 
and the additivity 
of every $\alpha$ is equivalent to the additivity of $\omega_v$. 

In Lemma~\ref{lem:a.11c} we have already seen that, 
for $v \in \cD_\g$, the continuity of the 
map $\omega_v \: \g \to V$ implies its linearity. 
Any more direct proof of the continuity of $\omega_v$ for 
$v \in \cD_\g$ would therefore lead to a more direct 
proof of Theorem~\ref{thm:additive}. 
\end{prob}

\begin{prob} For ``selfadjoint'' representations of the complex 
enveloping algebra $\cU(\g)_\C$ of a finite dimensional Lie algebra 
$\g$ on the dense subspace $\cD$ of the Hilbert space $\cH$, 
there exists an integrability criterion 
due to Powers (\cite[Thm.~4.5]{Po74}). The requirement is that 
the map $\pi \: \cU(\g)_\C \to \End(\cD)$ is ``completely strongly 
positive'' with respect to a certain convex cone $Q \subeq \cU(\g)_\C$. 
It would be very interesting to see if this result extends to 
Banach--Lie groups. 
\end{prob}

\begin{prob} As we have seen in Corollary~\ref{cor:7.4}, a 
positive definite function on a Lie group $G$ is smooth 
if it is smooth in some identity neighborhood. In some case 
one may even expect that the whole function can be reconstructed 
from the restriction to some identity neighborhood, even if it is 
not analytic. Theorem~2.1 in \cite{Jo91} contains a criterion for the 
extendability of a ``local'' positive definite function to the 
whole group for finite dimensional unimodular Lie groups. 
It would be very interesting to understand if there are variants of this 
result for more general topological groups and in particular for 
infinite dimensional Lie groups. Here the key point is to find 
appropriate positivity conditions, such as the 
{\it complete strong positivity} used in \cite[Cor.~4.1]{Jo91}. 
\end{prob}

{\bf Acknowledgment:} We thank St\'ephane Merigon and Christoph 
Zellner for a careful 
reading of this paper and for several discussions on its subject matter.

\end{document}